\documentclass[final]{siamltex}
\usepackage{times}
\usepackage{ucs} 
\usepackage[utf8x]{inputenc}
\usepackage{amsmath,amsfonts,amssymb,amsopn,amscd}
\usepackage{graphicx,color,colordvi}
\usepackage{hyperref}
\usepackage{slashbox}

\newtheorem{remark}{\it Remark\/}

\newcommand{\mK}{\mathbf{K}}
\newcommand{\mA}{\mathbf{A}}
\newcommand{\mB}{\mathcal{B}}
\newcommand{\mC}{\mathbf{C}}
\newcommand{\mF}{\mathcal{F}}
\newcommand{\mG}{\mathcal{G}}
\newcommand{\maH}{\mathcal{H}}
\newcommand{\mI}{\mathbf{I}}
\newcommand{\mJ}{\mathbf{J}}
\newcommand{\mn}{\mathbf{n}}
\newcommand{\mm}{\mathbf{m}}
\newcommand{\mr}{\mathbf{r}}
\newcommand{\mR}{\mathbf{R}}
\newcommand{\R}{\mathbb{R}}
\newcommand{\mcR}{\mathcal{R}}
\newcommand{\mU}{\mathbf{U}}
\newcommand{\mV}{\mathbf{V}}
\newcommand{\mv}{\mathbf{v}}
\newcommand{\mW}{\mathbf{W}}

\newcommand{\mL}{\mathbf{L}}
\newcommand{\mS}{\mathbf{S}}
\newcommand{\mX}{\mathbf{X}}
\newcommand{\suml}{\sum\limits}
\newcommand{\vf}{\mV^f_\lambda}
\newcommand{\norm}[1]{\left \| #1 \right \|}

\usepackage{ifthen}
\newboolean{DisplaySOS}
\setboolean{DisplaySOS}{true} 
\newcommand{\SOS}[1]{\ifthenelse{\boolean{DisplaySOS}}{{\bf[#1]}}{}}
\title{Local/global analysis of the stationary solutions of some neural field equations\thanks{This work was partially supported by the ERC grant 227747 - NERVI}}
\author{Romain Veltz\thanks{IMAGINE/LIGM, Universit\'e Paris Est and NeuroMathComp team.} \and Olivier Faugeras\thanks{Corresponding author, NeuroMathComp team, INRIA, CNRS, ENS Paris, olivier.faugeras@sophia.inria.fr}}

\begin{document}

\maketitle


\pagestyle{myheadings}
\thispagestyle{plain}
\markboth{R. VELTZ AND O. FAUGERAS}{STATIONARY SOLUTIONS OF NEURAL FIELDS EQUATIONS}


\begin{abstract}
Neural or cortical fields are continuous assemblies of mesoscopic models, also called neural masses, of neural populations
that are fundamental in the modeling of macroscopic parts of the brain. 
Neural fields are described by nonlinear integro-differential equations. 
The solutions of these equations represent the state of activity of these populations when submitted to  inputs from
neighbouring brain areas. Understanding the properties of these solutions is essential in advancing our understanding of
the brain.
In this paper we study the dependency of the stationary
solutions of the neural fields equations with respect to the stiffness of the nonlinearity and the contrast of the external inputs. This is done
by using degree theory and bifurcation theory in the context of functional, in particular infinite dimensional, spaces.
The joint use of these two theories allows us to make new detailed predictions about the global and local behaviours of the solutions.
We also provide a generic finite dimensional approximation of these equations which allows
us to study in great details two models. 
The first model is a neural mass model of a cortical hypercolumn of orientation sensitive neurons, the ring model \cite{shriki-hansel-etal:03}. The second model is a general neural field model where the spatial connectivity is
described by heterogeneous Gaussian-like functions. 
\end{abstract}
\section{Introduction}
Neural or cortical fields are continuous assemblies of mesoscopic models, also called neural masses, of neural populations
that are essential in the modeling of macroscopic parts of the brain. 

They were first studied by Wilson and Cowan \cite{wilson-cowan:73,amari:77} and
are fundamental in the design of the models of the visual cortex proposed by Bressloff (\cite{bressloff-cowan-etal:01}).
Neural fields describe the mean activity of neural populations which are described by nonlinear integro-differential equations. 
The solutions of these equations represent the state of activity of these populations when submitted to  inputs from
neighbouring brain areas. Understanding the properties of these solutions is important for advancing our understanding of
the brain.

Among these solutions, the persistent states (or stationary solutions) are interesting for at least two reasons. First, when dealing with autonomous systems, looking for persistent states helps to understand the dynamics because they are an easy way to divide the phase space into smaller components. Moreover, as the connectivity is often chosen symmetric (see the review by \cite{ermentrout:98}), the dynamics is described by heteroclinic orbits. Second, they are thought to model well the memory holding tasks on the time scale of the second which has been demonstrated by experimentalists on primates \cite{colby-duhamel-etal:95,funahashi-bruce-etal:89,miller-erickson-etal:96}.

There are several items that arise in the mathematical description of these equations. The domain $\Omega$ of integration (typically a piece of cortex) which is of dimension $d$ and can be bounded or unbounded, the type of the nonlinearity (sigmoidal, Heaviside), the type of the connectivity function that appears in the integral (homogeneous, i.e. translation invariant, or heterogeneous), and the number $p$ of neuronal populations that are modeled.

The nonlinearity in the neural field equations is most of the time a Heaviside function which leads to several mathematical difficulties,  one of them being being to specify the correct functional space. We found easier in \cite{faugeras-veltz-etal:09} to use smooth nonlinear functions instead, to study the persistent states. Only a few papers use sigmoidal functions, e.g., \cite{atay-hutt:05,venkov-coombes-etal:07},  or an alpha function\footnote{$e^{-1/v^2}H(v)$, $H$ the Heaviside function.}, \cite{laing-troy:03}. 
In this paper we make the assumption that the nonlinearity is a sigmoid function, i.e. infinitely differentiable. Also, in this paper, the only assumption on the connectivity kernel is that it is square integrable, hence it can be heterogeneous (not translation invariant) as it is often assumed.

Surprisingly, there are few papers dealing with the persistent states. Their authors use two main methods: Turing patterns (or bifurcation theory) and reduction to ODEs and PDEs.
Bifurcation theory helps to understand the local structure of the set of persistent states and also gives the local dynamical structure. 
The first method is used in almost every paper, e.g. \cite{atay-hutt:05,coombes-owen:04,blomquist-wyller-etal:05,bressloff-folias-etal:03,venkov-coombes-etal:07}, using a translation invariant connectivity function (hence a convolution) and leads to the descrition of a lot of behaviours depending on relations involving the spatial frequency $k$ and the time frequency $\omega$: traveling waves, breathers, persistent states\ldots The numerical computation of these states depends on their stability.  We generalize this approach as follows: we do not care about the spatial structure of the cortical states (indexed by $k$ in previous studies) but consider a cortical state $\mV(\cdot)$ as a point of a (functional) vector space, and how this point varies (bifurcates) when the relevant parameters in the neural field equations vary. By doing so we are able to harvest a lot more results that do not depend upon the translation invariant assumption but of course also apply to this case.

The second method is to reduce the persistent state equation to the problem of finding a homoclinic orbit to ODEs \cite{laing-troy-etal:02} when $\Omega=\mathbb R$ and $p=1$. If $\Omega=\mathbb R^2$, it reduces to the finding of homoclinic orbits to PDEs \cite{laing-troy:03}. The obvious advantage is that they can use finite-dimensional tools \cite{kuznetsov:98,guckenheimer-holmes:83} or PDEs methods, e.g. \cite{kielhofer:03}, for the bifurcation analysis. Hence, the authors are able to compute the persistent states independently of their stability and to show their number as a function of the strength of the connections.

Working with an unbounded domain $\Omega$ is not biologically relevant and also raises some mathematical difficulties. We work with a bounded $\Omega$.
Besides  its biological relevance, this hypothesis is \textbf{crucial} for our  mathematical analysis: it implies that there is at least one persistent state for any set of parameters and also it provides bounds for these persistent states.

In a previous paper \cite{faugeras-veltz-etal:09}, we started the analysis of the neural fields equations defined over a finite part of the cortex from two different viewpoints, theoretical and numerical. We proved some results concerning the dynamics and gave a method to efficiently compute the persistent state (or stationary solution) under the assumption that it was unique.

In the same article, we unsurprisingly found that the dynamics was poor (every initial condition converged to a single persistent state) if the stiffness of the nonlinearity was small, for example the system could not exhibit oscillatory behaviours. More importantly we did not give a way to compute the persistent states when the system featured more than one of such state.

In this paper we relax the hypothesis on the uniqueness of the persistent state and try to understand the structure of the set, noted $\mathcal B$, of stationary solutions as the parameters\footnote{We choose the stiffness of the nonlinearity and the the contrast of the external inputs but our analysis applies to other parameters.} vary over large scales, hence we are not only interested in local behaviours. The local structure of the set $\mathcal B$ is understood using bifuraction theory (see \cite{ma-wang:05,kielhofer:03,haragus-iooss:09}) whereas  degree theory helps to understand its global structure. This latter theory predicts stationary solutions that cannot be predicted using only bifurcation theory. Notice that bifurcation theory also gives the local dynamics. In order to compute numerically these persistent states, we propose a multiparameter continuation scheme that allows to compute non-connected branches\footnote{A branch is a one-dimensional set of stationary solutions obtained by varying one parameter.} of persistent states
The case when there is no inputs can be (at least locally) done analytically and is used to predict and verify the outcome of the numerical computations.

For conducting our numerical experiments we consider a very general class of approximating connectivity kernels, the Pincherle-Goursat kernels \cite{tricomi:85}, which reduce exactly the dynamics to a system of ODEs whose dimension is directly related to the level of approximation and can be arbitrarily large. Hence our choice to use infinite dimensional techniques to the integral equation is guided by the fact that it offers a simple, albeit abstract, conceptual framework in which the behaviours of interest to us can be described in a clean and dimension-independent manner. When it comes to numerical experiments we use  the system of ODEs provided by the Pincherle-Goursat kernels. 

The paper is organized as follows. In section \ref{section:general} we introduce a very general functional framework for studying the neural field equations which allows us to pose the problem as a Cauchy problem in a functional space and to derive a number of useful properties for the (stationary) solutions of these equations. In section \ref{section:bifurcations} we combine the results of the previous section with a bifurcation study in order to obtain more information about $\mathcal B$. A numerical scheme is proposed to compute the structure of $\mathcal B$. In section \ref{section:PG} we show how to reduce the neural field equations to a set of ODEs with an arbitrary precision through the use of the Pincherle-Goursat kernels. In section \ref{section:ring} we study in detail a neural mass model that reduces exactly to a finite set of ODEs, the ring model. In section \ref{section:numerics2d}, we compute the stationary solutions for a model of 2 populations of neurons on a bi-dimensional cortex connected by heterogenous Gaussian-like functions: it allows us to provide examples of the predictions obtained in sections \ref{section:general} and \ref{section:bifurcations}.

\section{General framework}\label{section:general}

We consider the following formal neural field equation defined over a bounded piece of cortex and/or feature space $\Omega\subset \mathbb R^d$. We wish to encompass in the same formalism the geometry of the cortex, seen as a bounded piece of $\R^d$, $d=1,2,3$, and the geometry of such feature spaces as edge or motion orientations (directions) which are represented by a value between 0 and $\pi$ ($2\pi$), or scale which is represented by a positive number, e.g. the spatial frequency. In section \ref{section:ring} we analyze in detail an example where the focus is on edge directions while in section \ref{section:numerics2d} we analyze another example in which the focus is on the two-dimensional geometry of the cortex.
\begin{equation}
\label{eq:NM}
\left\{
\begin{array}{lcl}
\dot\mV(\mr,t)&=&-\mL \cdot \mV(\mr,t)+\left[\mJ(t) \cdot \mS(\boldsymbol{\lambda} (\mV(t)-\boldsymbol{\theta}))\right](\mr)+\mI_{ext}(\mr,t) \quad t>0\\
\mV(\cdot,0)&=&\mV_0(\cdot)
\end{array}
\right.
\end{equation}
 This equation is an initial value problem that describes the time variation of the $p$-dimensional vector function $\mV$ defined on $\Omega$, starting from the initial condition $\mV_0$, a function defined on $\Omega$. At each time $t \geq 0$ $\mV$ belongs to some functional space, in effect a Hilbert space $\mF$, that we describe in the next section. We now discuss the various quantities that appear in \eqref{eq:NM}.

\noindent
$\mJ(t)$ is a linear operator from $\mF$ to itself defined by:
\begin{equation}\label{eq:W}
 [\mJ(t) \cdot \mV(t)](\mr)=\int_\Omega \mJ(\mr,\mr',t)\mV(\mr',t)\,d\mr',
\end{equation}
where $\mJ(\mr,\mr',t)$ is a $p \times p$ matrix that describes the ``strength'' of the connections. We also describe in the next section the functional space to wich $\mJ(t)$ belongs and the conditions it must satisfy in order for the equation \eqref{eq:NM} to be well-defined.

The external current input, $\mI_{ext}(\cdot,t)$,  is in $\mF$ for all $t \geq 0$.

The function 
$\mS : \R^p \to \R^p$ is defined by $\mS(x)=[S(x_1),\cdots,S(x_p)]^T$, where $S: \R \to (0,1)$ is the normalized sigmoid function of equation
\begin{equation}\label{eq:sig}
S(z)=\frac{1}{1+e^{-z}}.
\end{equation}
It is infinitely differentiable on $\R^p$ and all its derivatives $S^{(q)}(x)$, $q=1,2,\cdots$ are bounded. For all integer $q\geq 1$ we note $\mS^{(q)}(x)$ the $p \times p$ diagonal matrix ${\rm diag}(S^{(q)}(x_1),\cdots,S^{(q)}(x_p))$. Because of the form of the function $\mS$, the $q$th order derivative of $\mS$ at $x \in \R^p$ is the multilinear function defined by
\begin{equation}\label{eq:Sderiv}
 D^q\mS(x) \cdot (y_1,\cdots,y_q)=\mS^{(q)}(x) \cdot (y_1 \cdots y_q),\quad y_i \in \R^p,i=1,\cdots,q
\end{equation}
where $y_1 \cdots y_q$ is the component pointwise product of the $q$ vectors $y_1,\cdots,y_q$ of $\R^p$, i.e. the vector of $\R^p$ whose $k$th coordinate, $k=1,\cdots,p$ is equal to the product of the $q$ $k$th coordinates of each vector $y_i$, $i=1,\cdots,q$.

$\boldsymbol{\lambda}$ is the $p \times p$ diagonal matrix ${\rm diag}(\lambda_1,\cdots,\lambda_p)$, $\lambda_i \geq 0$, $i=1,\cdots,p$ that determines the slope of each of the $p$ sigmoids at the origin. We note $\lambda_m$ the
maximum value of the $\lambda_i$s.

$\boldsymbol{\theta}$ is a $p$-dimensional vector that determines
the threshold of each of the $p$ sigmoids, i.e. the value of the membrane potential corresponding to $50\%$ of the maximal activity.

The diagonal $p \times p$ matrix $\mL$ is equal to ${\rm diag}(\frac{1}{\tau_1},\cdots,\frac{1}{\tau_p})$, where the positive numbers $\tau_i$, $i=1,\cdots,p$  determine the exponential decrease dynamics of each neural population.

As recalled in the introduction and detailed in \cite{faugeras-grimbert-etal:08,faugeras-touboul-etal:09} this equation corresponds to a mesoscopic description of population of neurons which is called voltage-based by Ermentrout \cite{ermentrout:98}. $\mV$ in equation \eqref{eq:NM} therefore has the biological interpretation of an average membrane potential of $p$ populations of neurons. As pointed out by the same authors, there is also an activity-based mesoscopic description which leads to the following initial value problem:
\begin{equation}
\label{eq:NMa}
\left\{
\begin{array}{lcl}
\dot\mA(\mr,t)&=&-\mL_a \cdot \mA(\mr,t)+\mS\left(\boldsymbol{\lambda}\left(\left[\mJ(t) \cdot \mA\right](\mr,t)+\mI_{ext}(\mr,t)\right)\right) \quad t>0\\
\mA(\cdot,0)&=&\mA_0(\cdot)
\end{array}
\right.
\end{equation}
$\mL_a\neq \mL$ (see \cite{ermentrout:98}) because they do not have the same biologiocal meaning: One is related to the synaptic time constant and the other to the cell membrane time constant. We let $\mL_a={\rm diag}(\alpha_1,\cdots,\alpha_p)$.

The two problems \eqref{eq:NM} and \eqref{eq:NMa} are closely related. In particular there is a one to one correspondence between their equilibria, as recalled below for the ring model of section \ref{section:ring}, which is an activity-based model.
\subsection{Choice of the appropriate functional space}\label{subsection:functional}
The problem at hand is to find an appropriate mathematical setting, i.e. to choose the functional space $\mF$, for the neural field equations based on three criteria 1) problems \eqref{eq:NM} and \eqref{eq:NMa} should be well-posed, 2) its biological relevance and 3) its ability to allow  numerical computations. 

A Hilbert space is appealing because of its metric structure induced by its inner product. Hence, a natural choice arising from the definition of the linear operator $\mJ$ would be the space  ${\rm \mL}^2(\Omega,\mathbb R^p)$. 
But this space allows the average membrane potential to be singular which is biologically excluded. For example the function $\mr\to|\mr|^{-1/2}\in{\rm L}^2(\Omega,\mathbb R^p)$ if $d\geq 0$. It would be desirable that this potential  be bounded on the cortex. A way to achieve this is to allow for more spatial regularity of the membrane potential: it is \textbf{reasonable}, E.g., from optical imaging measurements. to choose $\mr\to\mV(\mr)$ as being differentiable almost everywhere.

For technical reasons that will become clear later we choose the Sobolev space 
\[
\mF=\mW^{m,2}(\Omega,\mathbb R^p), \quad m \in \mathbb{N},
\]
with the inner product:
\begin{equation}
 \label{eq:psF}
\left\langle \mX_1,\mX_2\right\rangle_{\mF}=\sum_{|\alpha|=0}^m\left\langle D^\alpha\mX_1,D^\alpha\mX_2\right\rangle_{{\mL}^2(\Omega,\mathbb R^p)},\ \forall \mX_1,\,\mX_2 \in \mF,
\end{equation}

where, as usual, the multi-index $\alpha$ is a sequence $\alpha=(\alpha_1,\cdots,\alpha_d)$ of positive integers, 
$
 |\alpha|=\sum_{i=1}^d \alpha_i,
$
and the symbol $D^\alpha$ represents a partial derivative:
\[
 D^\alpha \varphi=\frac{\partial^{\alpha_1+\cdots+\alpha_d}}{\partial r_1^{\alpha_1}\cdots \partial r_d^{\alpha_2}}\varphi,
\]
where $\varphi: \Omega \to \R$ and $\mr=(r_1,\cdots,r_d)$.

Note that $\mF=\underset{p}{\underbrace{W^{m,2}(\Omega) \times \cdots \times W^{m,2}(\Omega)}}$. Later we use the notation $\mG$ for $W^{m,2}(\Omega)$, i.e. $\mF=\mG^p$.
\subsubsection{Choosing the value of $m$}

The value of the integer $m$ determines the regularity of the functions $\mV$ that represent the cortical states as well as that of the connectivity function $\mJ$. It turns out that, depending upon the relative values of $m$ and the dimension $d=1,2,3$ of the space where we represent the cortical patch $\Omega$ and/or the feature space, $\mF$ is a (commutative) Banach algebra for the pointwise multiplication. This is important when using bifurcation theory because we need to use Taylor expansions of the right-hand-side of (\ref{eq:NM}).

Being a Banach algebra property requires some regularity of $\Omega$, i.e. of its boundary $\partial \Omega$. Because in the numerical experiments of section \ref{section:numerics2d} we work with the open square $\Omega=[-1,\ 1]^2$, we will assume the weakest regularity, the so-called cone property, see \cite[Chapter IV]{adams:75} for a definition. In practice one can safely assume that $\partial \Omega$ is $C^1$. Roughly speaking it means that this boundary is a $C^1$-manifold of dimension $d-1$ of $\R^d$. The exact definition, called the uniform $C^1$-regularity, can again be found in \cite[Chapter IV]{adams:75}. It implies some weaker regularity assumptions, such as the cone property.

Under this assumption on $\Omega$ we can adapt the theorem in \cite[Chapter V, Theorem 5.24]{adams:75} to obtain the following result :

\begin{corollary}\label{corollary:banach}
If the relation $2m > d$ holds, then
 $\mF$ is a commutative Banach algebra with respect to component pointwise multiplication, i.e. for all $\mU$ and $\mV$ elements of $\mF$ their component pointwise product $\mU\mV$ is in $\mF$ and there exists a positive constant
$K^*$ such that
\[
 \norm{\mU\mV}_\mF \leq K^* \norm{\mU}_\mF\,\norm{\mV}_\mF \quad \forall \mU,\,\mV \in \mF,
\]
where $\mU\mV=(\mU_1\mV_1,\cdots,\mU_p\mV_p)$.
\end{corollary}

Let us explore some consequences of this proposition on our possible choices for the value of $m$. These choices are guided by the idea of constraining the functional space $\mF$ as little as possible, given the fact that $\mW^{m,2} \subset \mW^{n,2}$ for all integers $0 \leq n < m$
\begin{description}
 	\item[$d=1$:] The relation $2m > 1$ implies $m \geq 1$ and we choose $m=1$.
	\item[$d=2,3$:] The relation $2m > d$ implies $m \geq 2$ and we choose $m=2$. 
\end{description}
Similar results hold for higher values of the dimension $d$ but these cover the cases discussed in this article.

We summarize all this in the following proposition
\begin{proposition}[Choice of $\mF$]\label{proposition:mF}
If $\Omega$ is has the cone property, in particular if it is uniformly $C^1$-regular, then we have $\mF=\mW^{1,2}(\Omega,\R^p)$ for $d=1$ and $\mF=\mW^{2,2}(\Omega,\R^p)$ for $d=2,3$. In all three cases, $\mF$ is a commutative Banach algebra with respect to component pointwise multiplication.
\end{proposition}

\subsubsection{The choice of $\mJ$}

We assume that $\mJ(\cdot,\cdot,t) \in  {\rm \mW}^{m,2}(\Omega \times \Omega,\mathbb R^{p \times p})$ for all $t>0$.

As a consequence, the Frobenius $\mF$-norm of the linear operator $\mW(t)$, noted $\norm{\mJ(t)}_\mF$, is well-defined
\[
\norm{\mJ(t)}_\mF^2=\sum_{|\alpha|,\,|\alpha'|=0}^m\int_{\Omega \times \Omega} \norm{D^\alpha D^{\alpha'}\mJ(\mr,\mr',t)}_F^2\,d\mr\,d\mr',
\]
where it is understood that the partial derivative operator $D^\alpha$ (respectively $D^{\alpha'}$) acts on the variable $\mr$ (respectively $\mr'$) and $\norm{\mJ(\mr,\mr',t)}_F$ is the Frobenius norm of the matrix $\mJ(\mr,\mr',t)$,
\[
\norm{\mJ(\mr,\mr',t)}_F^2=\sum_{i,j} J_{ij}(\mr,\mr',t)^2
\]

This hypothesis also ensures the existence of the operator $\mJ$ as a linear operator from $\mF$ to $\mF$.
\begin{proposition}\label{prop:FtoF}
 Assume that $\mJ(\cdot,\cdot,t) \in  {\rm \mW}^{m,2}(\Omega \times \Omega,\mathbb R^{p \times p})$ for all $t>0$. Then
equation \eqref{eq:W} defines a linear operator from $\mF$ to $\mF$.
\end{proposition}
\begin{proof}
see appendix \ref{appendix:Banach}
\end{proof}

\

\noindent
Since $\mS$ is a $\R^p$-valued bounded function on $\Omega$ and such functions are in $\mL^2(\Omega)$ the right-hand side of \eqref{eq:NM} is an element  of $\mF$. Similarly, because $\mS$ is also infinitely differentiable and all its derivatives bounded, the right-hand side of \eqref{eq:NMa} is an element  of $\mF$.

Finally we  note $\norm{|\mJ(t)|}$ the operator norm of $\mJ(t)$, i.e.
\[
 \sup_{\norm{\mV}_\mF \leq 1} \frac{\norm{\mJ(t) \cdot \mV}_\mF}{\norm{\mV}_\mF}
\]
It is known, see e.g. \cite{kato:95}, that
\[
 \norm{|\mJ(t)|} \leq \norm{\mJ(t)}_\mF
\]

\begin{remark}
Note that the relation $2m > d$ and our regularity assumption on $\Omega$, imply, through the so-called Sobolev imbedding theorem, \cite[Chapter V, lemma 5.15]{adams:75}, that 
\[
	\mG \rightarrow C^0_B(\Omega),
\]
where $C^0_B(\Omega)$ is the set of continuous bounded functions on $\Omega$.
The consequence is that the state vectors $\mV$ and the connectivity matrix $\mJ$ are essentially bounded, which has certainly the right biological flavor.
\end{remark}

\subsection{The Cauchy problem}
In this section we prove that equation \eqref{eq:NM} is well-posed and provide some properties of its solutions.

We rewrite it as a Cauchy Problem, i.e. as an ordinary differential equation on the functional space $\mF$. This turns out to be convenient for the upcoming computations.
\begin{equation}
\label{eq:ChauchyPb}
\left\{
\begin{array}{lcl}
\frac{d\mV}{dt}&=&-\mL \cdot \mV+\mR(t,\mV) \quad t>0\\
\mV(0)&=&\mV_0\in\mF
\end{array}
\right.
\end{equation}
The nonlinear operator $\mR$ is defined by 
\begin{equation}\label{eq:R}
 \mR(t,\mV)=\mJ(t) \cdot \mS(\boldsymbol{\lambda} (\mV-\boldsymbol{\theta}))+\mI_{ext}(t)
\end{equation}

Proposition \ref{prop:FtoF} shows that $\mR(t,\cdot):\mF\to\mF$ for all $t>0$. We  have the further properties:
 \begin{lemma} \label{lemma:R}
$\mR$ satisfies the following properties :
\begin{itemize}

               \item $\forall\ \text{integer} \ q\geq 0,\ \ \mR(t,\cdot)\in C^q(\mF,\mF)$
  and ${D^q\mR(t,\mV_0)}=\boldsymbol{\lambda}^q{\mJ(t) \mS^{(q)}(\boldsymbol{\lambda}(\mV_0-\boldsymbol{\theta}))}$ for all $\mV_0$ in $\mF$.
               \item $\norm{\mR(t,\mU_1)-\mR(t,\mU_2)}_\mF\leq\lambda_m\norm{\mJ(t)}_\mF\norm{\mU_1-\mU_2}_\mF$ for all $t>0$ and for all $\mU_1$, $\mU_2$ in $\mF$.

               \item $\mR(t,\cdot)$ is a compact operator on $\mF$ for all $t>0$.
              \end{itemize} 
 \end{lemma}
 \begin{proof}
 It is easy to see from the definition \eqref{eq:R} of $\mR$ that, if it exists, $D^q\mR(t,\mV_0)[\mU_1,\cdots,\mU_q]=\boldsymbol{\lambda}^q\mJ(t) \cdot \left(\mS^{(q)}(\boldsymbol{\lambda}(\mV_0-\boldsymbol{\theta}))\cdot (\mU_1 \cdots \mU_q)\right)$. The notation $\mU_1 \cdots \mU_q$ is the same as in the definition of $D^q\mS$ in equation \eqref{eq:Sderiv}, i.e. the component pointwise product of the $q$ functions $\mU_1,\cdots,\mU_q$ of $\mF$.

The $q$-multilinear operator $D^q\mR(t,\mV_0)$ is well-defined because, according to corollary \ref{corollary:banach}, $\mU_1 \cdots \mU_q$ is in $\mF$. The first property follows immediately.

It remains to show that the $q$-multilinear operator $D^q\mR(t,\mV_0)$ is continuous. We apply corollary \ref{corollary:banach} once more to show that 
\[                 
\norm{D^q\mR(t,\mV_0) \cdot (\mU_1 \cdots \mU_q)}_\mF \leq \norm{\boldsymbol{\lambda}^q\mJ(t)\mS^{(q)}(\boldsymbol{\lambda}(\mV_0-\boldsymbol{\theta})) }_\mF\norm{\mU_1 \cdots \mU_q}_\mF\leq K \prod_i\norm{\mU_i}_\mF
\]
for some positive constant $K$, that is $D^q\mR$ is continuous and $\mR(t,\cdot)\in C^q(\mF,\mF)$.

 The second property is proved in \cite{faugeras-veltz-etal:09}.

 For all $\alpha$, $|\alpha|\leq m$ $\mU\to\partial^\alpha\mJ \cdot \mU$ are compact operators on ${\rm \mL}^{2}(\Omega,\mathbb R^p)$ because they are Hilbert-Schmidt operators (see \cite[chapter X.2]{yosida:95}). Hence for any bounded $\mV_n$ sequence in $\mF$ hence in ${\rm \mL}^{2}(\Omega,\mathbb R^p)$, there exists a subsequence $\mV'_n$ such that $D^\alpha\mJ(t)\cdot\mS(\boldsymbol{\lambda}(\mV_n'-\boldsymbol{\theta}))$ are convergent in ${\rm \mL}^{2}(\Omega,\mathbb R^p)$ for $|\alpha|\leq m$. Then $\mR(t,\mV_n')$ is convergent in $\mF$ because of (\ref{eq:psF}). We have proved that $\mR(t,\cdot)$ is a compact operator on $\mF$.
 \end{proof}

\

Following closely \cite{faugeras-grimbert-etal:08} we have the following proposition.
\begin{proposition}
If the following two hypotheses are satisfied:
\begin{enumerate}
     \item The connectivity function $\mJ$ is in $C(\R^+; {\mW}^{m,2}(\Omega \times \Omega,\R^{p \times p})$ for the values of $m$ given in proposition \ref{proposition:mF}, and bounded, $\norm{\mJ(t)}_\mF\leq J$, $t\geq 0$,
     \item the external current $\mI_{\rm ext}$ is in $C(\R^+;\mF)$,
\end{enumerate}
then for any function $\mV_0$ in $\mF$ there is a unique solution $\mV$, defined on $\R^+$  and
continuously differentiable, of the  initial value problem \eqref{eq:NM}.
\end{proposition}
 This solution depends upon  $3p$ parameters, the slopes $\boldsymbol{\lambda}$, the thresholds $\boldsymbol{\theta}$ and the diagonal matrix $\mL$.

Even if we have made progress in the formulation of the Neural Field equations, there still remains the unsatisfactory possibility that the membrane potential  becomes unbounded as $t\to\infty$. However this is not the case as shown in the next proposition:

\begin{proposition}
\label{prop:bounded}
 If the external current is bounded in time $\norm{\mI_{\rm ext}}_\mF\leq I_{\rm ext}$, for all $t \geq 0$, then the solution of equation \eqref{eq:ChauchyPb} is bounded for each initial condition $\mV_0 \in \mF$.
\end{proposition}
\begin{proof}
Let us define $f: \R\times\mF \to \R^+$ as 
\[
f(t,\mV)\overset{\rm def}{=}\left\langle -\mL \cdot \mV+\mJ(t) \cdot \mS(\boldsymbol{\lambda}(\mV-\boldsymbol{\theta}))+\mI_{ext}(t),\mV\right\rangle_{\mF}=\frac{1}{2}\frac{d \| \mV \|_\mF^2}{d\,t}.
\]
We note $\tau_{max}=\max_{i=1\cdots,p} \tau_i$ and notice that
\[
f(t,\mV)\leq 
-\frac{1}{\tau_{max}}\|\mV\|_\mF^2+(J+I_{ext})\|\mV\|_\mF.
\]
Thus, if $\|\mV\|_\mF \geq 2\tau_{max}(W+I_{ext}) \overset{\rm def}{=} R$, $f(t,\mV)\leq-2 \tau_{max} (W+I_{ext})^2\overset{\rm def}{=}-\delta<0$.

Let us  show that the open ball of $\mF$ of center 0 and radius $R$, $B_R$, is stable under the dynamics of equation (\ref{eq:NM}). We know that $\mV(t)$ is defined for all $t \geq 0$s and that $f<0$ on $\partial B_R$, the boundary of $B_R$. 
We consider three cases for the initial condition $\mV_0$.

If $\mV_0\in B_R$ and set $\tau=\sup\left\lbrace t\,|\,\forall s\in\left[0,t\right], \mV(s)\in \overline{B}_R\right\rbrace$. Suppose that $\tau\in\mathbb R$, then $\mV(\tau)$ is defined and belongs to $\overline{B}_R$, the closure of $B_R$, because $\overline{B}_R$ is closed, in effect to $\partial B_R$. We also have $\frac{d}{dt}\|\mV\|_\mF^2|_{t=\tau}=f(\tau,\mV(\tau))\leq -\delta<0$ because $\mV(\tau)\in\partial B_R$. Thus we deduce that for $\varepsilon > 0$ and small enough, $\mV(\tau+\varepsilon)\in\overline{B}_R$ which contradicts the definition of $\tau$. Thus $\tau\notin\mathbb R$ and $\overline{B}_R$ is stable. 

Because $f<0$ on $\partial B_R$, $\mV_0\in\partial B_R$ implies that $\forall t > 0,\ \ \mV(t)\in B_R$.

Finally we consider the case $\mV_0\in \complement \overline{B}_R$. Suppose that $\forall t>0,\ \mV(t)\notin \bar B_R$, then $\forall t>0,\ \frac{d}{dt}\|\mV\|_\mF^2\leq-2\delta$, thus $\|\mV(t)\|_\mF$ is monotonically decreasing and reaches the value of $R$ in finite time when $\mV(t)$ reaches $\partial B_R$. This contradicts our assumption. Thus $\exists\tau >0 \,|\,\mV(\tau)\in B_R$.
\end{proof}
\noindent

\begin{corollary}
 If $\mV_0\notin B_R$ and $T=\inf\left\lbrace t>0|\mV(t)\notin B_R\right\rbrace $. Then $$T\leq \frac{\|\mV_0\|_\mF^2-R^2}{2\delta}$$
\end{corollary}

This proposition shows that $\overline{B}_R$ is an attracting set and that it suffices to study the dynamics within this set to have the long time behavior of the solutions of the Neural Fields Equations. 

This attracting set contains the stationary solutions of (\ref{eq:NM}), we devote the next section to their study. We quote a result from \cite{faugeras-veltz-etal:09} concerning their stability:
\begin{proposition}\label{proposition:Lyapunov}
 If the condition
\[
 \lambda_m \rho(\mJ_s)<1,
\]
holds, then every stationary solution of (\ref{eq:NM}) is globally asymptotically stable.
$\lambda_m=\max \lambda_i$, $\mJ_s$ is the symmetric part $(\mJ+\mJ^*)/2$ of the operator $\mJ$, and $\rho(\mJ_s)$ its spectral radius.
We define $\lambda_L$ to be $\rho(\mJ_s)^{-1}$.
\end{proposition}

\noindent
Similar results hold for the activity-based model \eqref{eq:NMa}.

\subsection{Properties of the set of persistent states }\label{subsection:blambda}
We look at the equilibrium states, noted $\mV^f_\lambda$, of (\ref{eq:NM}), when $\mI_{ext}$ and $\mJ$ do not depend upon the time. Our goal is to estimate their number and, if possible, to compute them numerically, for a given set of parameters. 

It is quite demanding to do it at a given point in the parameter space except in some very special cases\footnote{When the nonlinearity $\mJ \cdot \mS(\boldsymbol{\lambda} (\vf-\boldsymbol{\theta}))$ is small compared to the linear part $\mL \cdot \mV_\lambda^f$, we know there exists a unique solution and how to compute it efficiently. This was proved in \cite{faugeras-veltz-etal:09}.}.  We note that when  $\boldsymbol{\lambda}=0$ (or $\mJ=0$), the stationary equation is trivially solved. Hence, we can think of deforming this trivial solution to a solution when $\boldsymbol{\lambda}\neq0,\mJ\neq 0$. This raises a number of questions. Does such a ``manifold'' of solutions exist \textit{i.e.} can we link the trivial solution to a solution for any given set of parameters? If yes, are there any other solutions? How do these ``manifolds'' look globally?  These questions concern global properties of the set of solutions (existence of branches, existence of intersection points, connectedness\ldots) are difficult to answer. We provide some partial answers in the remaining of the article.

\

Before going deeper in the analysis, we need to simplify the parameter space. The equilibria, also called bumps, or persistent states , are stationary solutions (independent of time) of 
\[
0=-\mL \cdot \vf+\mJ \cdot \mS(\boldsymbol{\lambda} (\vf-\boldsymbol{\theta}))+\mI_{\rm ext}, 
\]
\noindent
We redefine $\mJ$ as $\mL^{-1} \mJ$, $\mV$ as $\mV-\boldsymbol{\theta}$ and $\mI_{ext}$ as $\mL^{-1}\cdot\mI_{ext}+\boldsymbol{\theta}$ and restrict our study to:
\begin{equation}
\label{eq:bump1}
0=-\vf+\mJ \cdot \mS(\boldsymbol{\lambda}\vf)+\mI_{\rm ext}
\end{equation}

Still, equation (\ref{eq:bump1}) contains many parameters such as the ones describing $\mJ$ and $\mI_{\rm ext}$, or the slopes $\boldsymbol{\lambda}$. Which parameters to choose for the continuation method: $\boldsymbol{\lambda}$ or $\mJ$ ? We decide to set $\mJ$ and control $\boldsymbol{\lambda}$ for two reasons :
\begin{itemize}
 \item the stationary solutions are bounded for $\boldsymbol{\lambda}\in\mathbb R_+^p$, see proposition \ref{prop:bumps}.1 which is not the case when $\|\mJ\|\to\infty$. This proves to be useful numerically.
 \item previous studies usually use a Heaviside nonlinearity which is formally equivalent to our nonlinearity when $\boldsymbol{\lambda}\text{'='}\infty$, varying $\boldsymbol{\lambda}$ can thus bridge the gap with previous approaches.
\end{itemize}
As a matter of fact the techniques we are about to expose are applicable to any set of parameters with minor modifications.
Hence, we now focus on the influence of the slopes $\boldsymbol{\lambda}$ on the solutions of (\ref{eq:NM}). We make the assumption that they are all are equal to $\lambda$, $\boldsymbol{\lambda} =\lambda {\rm Id}_p$, $\lambda \geq 0$.  where ${\rm Id}_p$ is the $p \times p$ identity matrix. The equation becomes
\begin{equation}
\label{eq:bump}
0=-\vf+\mJ \cdot \mS({\lambda}\vf)+\mI_{\rm ext}\overset{\rm def}{=}-F(\vf,\lambda)
\end{equation}
It is clear that when $\boldsymbol{\lambda}=0$, the stationary equation is trivially solved by 
\[
 \mV^f_0\overset{\rm def}{=} \mJ \cdot \mS(0)+\mI_{\rm ext}=\frac{1}{2} \mJ \cdot \mathbf{1}+\mI_{\rm ext},
\]
where $\mathbf{1}$ is the $p$-dimensional vector with all coordinates equal to 1.
Let $\mB_\lambda$ be the set of solutions of equation \eqref{eq:bump} for a given slope parameter $\lambda$ :
$$\mB_\lambda=\left\{\mV\,|\,F(\mV ,\lambda)=0 \right\}$$
We next provide some properties of the sets $\mB_\lambda$.

\begin{figure}[htbp]
	\centering
         \includegraphics[width=\textwidth]{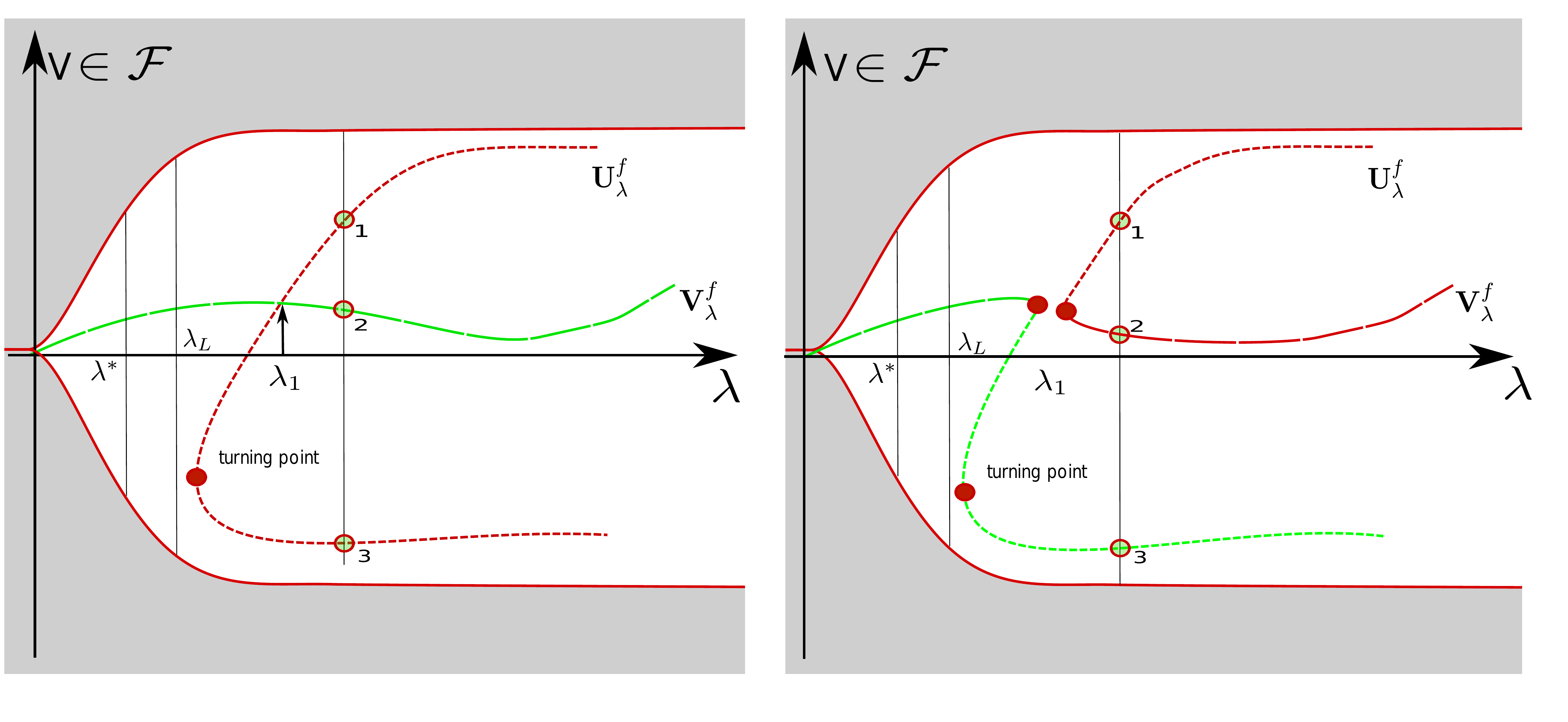}
       \caption{The white zone is the domain of existence of $\mV^f$. The grey zone is excluded thanks to proposition \ref{prop:bumps}.1. Left: non generic situation, a transcritical bifurcation occurs at $\lambda_1$. Right: same as left with a small perturbation (like a tiny change of the external input $\mI_{\rm ext}$), the transcritical bifurcation has opened up. $\lambda_L$ is defined in prop.\ref{proposition:Lyapunov} and $\lambda^*$ in prop.\ref{prop:bumps}.}
\label{fig:SigBorne}
\end{figure}

\begin{proposition} \label{prop:bumps}\ \\
\begin{enumerate}
 \item The persistent states satisfy the following inequality
\[
\left\|\mV_\lambda^f-\mV^f_0\right\|_\mF\leq \left\|\mJ\right\|_\mF \sqrt{p|\Omega|S_0\left(\frac{\lambda^2 B_1^2}{p|\Omega|}\right)},\text{ where } B_1\overset{\rm def}{=}\sqrt{p|\Omega|} \left\|\mJ\right\|_\mF+\norm{\mI_{\rm ext}}_\mF 
\]
where $S_0: \R \to \R$ is the ``shifted'' sigmoid defined by $S_0(x)=S(x)-S(0)$ and the constant $B_1$ is defined in proposition \ref{prop:ineq1} of appendix \ref{appendix:lemmas}.
\item If the condition
\begin{equation}\label{eq:unique}
\lambda \, \|\mJ\|_\mF<1
\end{equation}
is satisfied, then $\#\mathcal B_\lambda=1$. We define $\lambda^*$ to be $\norm{\mJ}_\mF^{-1}$.
\item $\forall\lambda\in\mathbb R^+,\ \mB_\lambda\neq\emptyset,$
 \item If we know an even number of solutions in $\mathcal B_\lambda$, then $\mathcal B_\lambda$ contains at least one more solution.
\item Let $0\leq a < b$ be two reals,  and consider the set $\mathcal{B}= \cup_{\lambda \in [a,\,b]} \mathcal{B}_\lambda \times \{\lambda\}$. Then $\mathcal{B}$ contains a connected component $\mathcal{C}$ which intersects $\mathcal{B}_{a} \times \{a\}$ and $\mathcal{B}_{b} \times \{b\}$.
\end{enumerate}
\end{proposition}
\begin{proof}
\begin{enumerate}
\item From lemma \ref{lemma:Sromain} in appendix  \ref{appendix:lemmas} we have $S_0(\lambda V^f_{\lambda i})^2 \leq S_0(\lambda^2 (V^f_{\lambda i})^2)$, $i=1,\cdots,p$. Therefore
\[
 \sum_{i=1}^p S_0(\lambda V^f_{\lambda i})^2 \leq \sum_{i=1}^p S_0(\lambda^2 (V^f_{\lambda i})^2) \leq p S_0\left( \frac{\lambda^2}{p}\sum_{i=1}^p (V^f_{\lambda i})^2 \right).
\]
The second inequality comes from Jensen's and the fact that $S_0(\cdot)$ is concave in $\R^+$. It then follows, using again Jensen's inequality and the fact that $\mS_0$ is monotonously increasing, that
\[
 \norm{\mS_0(\lambda \vf)}^2_{{\rm L}^2(\Omega,\mathbb R^p)} \leq p |\Omega| S_0\left(\frac{\lambda^2}{p|\Omega|} \norm{\vf}^2_{{\rm L}^2(\Omega,\mathbb R^p)}\right)\leq p |\Omega| S_0\left(\frac{\lambda^2}{p|\Omega|} \norm{\vf}^2_\mF\right)
\]
Now 
\begin{multline}
  \norm{\mV_\lambda^f-\mV^f_0}^2_\mF=\norm{\mJ\cdot\mS_0(\lambda \vf)}^2_\mF=\sum_{|\alpha|=0}^m\norm{D^\alpha\mJ\cdot\mS_0(\lambda \vf)}^2_{{\mL}^2(\Omega,\mathbb R^p)}\\
\leq \norm{\mS_0(\lambda \vf)}^2_{{\mL}^2(\Omega,\mathbb R^p)} \sum_{|\alpha|=0}^m\norm{D^\alpha\mJ}_{{\rm L}^2(\Omega,\mathbb R^p)}
\leq \norm{\mS_0(\lambda \vf)}^2_{{\rm L}^2(\Omega,\mathbb R^p)} \norm{ \mJ}_\mF^2
\end{multline}
The inequality then follows from proposition \ref{prop:ineq1}.
\item Use the Picard Theorem. As shown in figure \ref{fig:SigBorne}, this imposes that $\lambda^* \leq \lambda_L$. Indeed,
as $\rho(\mJ_s)\leq \||\mJ_s\||$ and $\norm{|\mJ_s|}=\norm{|\mJ|}$, and $\norm{|\mJ|} \leq \norm{\mJ}_\mF$, we have $\lambda^*=\norm{\mJ}_\mF^{-1}\leq \lambda_L=\rho(\mJ_s)^{-1}$.
\item The first property is that it is non empty:
 in \cite{faugeras-veltz-etal:09} we proved that persistent always exist in ${\mL}^2(\Omega,\mathbb R^p)$ for any positive values of $\lambda$. However, our ambient functional space is different since we  require more space regularity. Let us consider a persistent state $\vf\in {\mL}^2(\Omega,\mathbb R^p)$ given by \cite{faugeras-veltz-etal:09}. It satisfies $\vf=\mR(\vf)$. But $\mR:{\mL}^2(\Omega,\mathbb R^p)\to\mF$ because of the assumptions on $\mJ(\mr,\mr')$ and $\mI_{ext}$. Hence any persistent state in ${\mL}^2(\Omega,\mathbb R^p)$ in fact belongs to $\mF$.
 \item Suppose that $\mathcal B_\lambda$ has an infinite number of solutions, then the proposition holds. If now, $\mathcal B_\lambda$ has a finite number of solutions, we can assume that these points are non critical. Then according to the Leray-Schauder degree theory sketched in appendix \ref{appendix:LS} we have
\begin{multline*}
{\rm deg}_{\rm LS}(F(\mV,\lambda),B_r,0)=\suml_{\mV^f_\lambda\in\mathcal B_\lambda}{\rm sign}\ {\rm det}_{\rm LS}(D_V F(\mV_\lambda^f,\lambda),B_r,0)=\\
\suml_{\mV^f_\lambda\in\mathcal B_\lambda}{\rm sign}\ {\rm det}_{\rm LS}({\rm Id}-\lambda \mJ D\mS(\lambda\mV^f_\lambda)),
\end{multline*}
where $r=2B_1$ and $B_1$ is defined in proposition \ref{prop:ineq1} in appendix  \ref{appendix:lemmas}.
We prove in  corollary \ref{coro:bounded} in appendix \ref{appendix:LS} that the first term is equal to 1.
Suppose now that $\mathcal B_\lambda$ contains an even number of points, say $2k$ among which $l$ correspond to a negative sign and hence $2k-l$ to a positive sign. The sum that appears in the last term is equal to $2k-2l$, hence even. Hence $\mathcal B_\lambda$ must possess an odd number of points. 
\item The proof uses the Leray-Schauder theorem, see appendix \ref{appendix:LS}. We apply the theorem to the function $F: \mF \times J \to \mF$ which is of the form ${\rm Id}+m$, with $m(\cdot)=-\mR$. Because $m$ is compact on $\mF \times J$ (see proof in \cite{faugeras-grimbert-etal:08}), $\mathcal{B}_{a} \times \{a\}$ bounded, and there exists an open  bounded neighbourhood $\mathcal{U}_{a}$ of $\mathcal{B}_{a}$ such that ${\rm deg}_{\rm LS}(F(\cdot,a),\mathcal{U}_{a},0) \neq 0$ (corollary \ref{coro:bounded} in appendix \ref{appendix:LS}), the conclusion follows since the connected component cannot be unbounded, $\mathcal{B}$ being bounded.
\end{enumerate}
\end{proof}

\

This proposition answers some of the previous questions. For any positive value of $\lambda$, there is always at least one persistent state and we can find a way to connect the trivial solution $\mV^f_0$ to a persistent state corresponding to an arbitrary value of the  parameter $\lambda$. We return to this connection later. All this is true if we choose, say $\|\mJ\|_\mF$, as a parameter instead of $\lambda$. We will also see that not all the solutions in $\mathcal B_\lambda$ are in the connected component of $(\mV_0^f,0)$ in $\mF \times \mathbb R_+$.

\

Regarding the connection between $\mV_0^f$ and $\mV$, proposition \ref{prop:bumps} does not give us any indication on its regularity but we have the following corollary.
\begin{corollary}
\label{coro:connected}
Let $a$ and $b$ be as in proposition \ref{prop:bumps}. For all $\varepsilon > 0$ there exists a finite sequence $(\mV_1,\lambda_1)$, \ldots, $(\mV_n,\lambda_n)$ of points of $\mathcal{C}$ such that $\norm{\mV_i-\mV_{i+1}}_\mF \leq \varepsilon$ for $i=1,\cdots,n-1$ and $\lambda_1=a$, $\lambda_n=b$.
\end{corollary}
\begin{proof}
$\mathcal{C}$ is connected for any topology equivalent to the product topology of $\mF \times [a,\,b]$, e.g. for the metric defined by $d((\mV_1,\lambda_1),(\mV_2,\lambda_2))=\norm{\mV_1-\mV_2}_\mF+|\lambda_1-\lambda_2|$. Since it is connected for this metric, it is also well-chained \cite{choquet:69}, and the conclusion follows.
\end{proof}

\
 
\begin{figure}[htbp]
\centering
         \includegraphics[width=\textwidth]{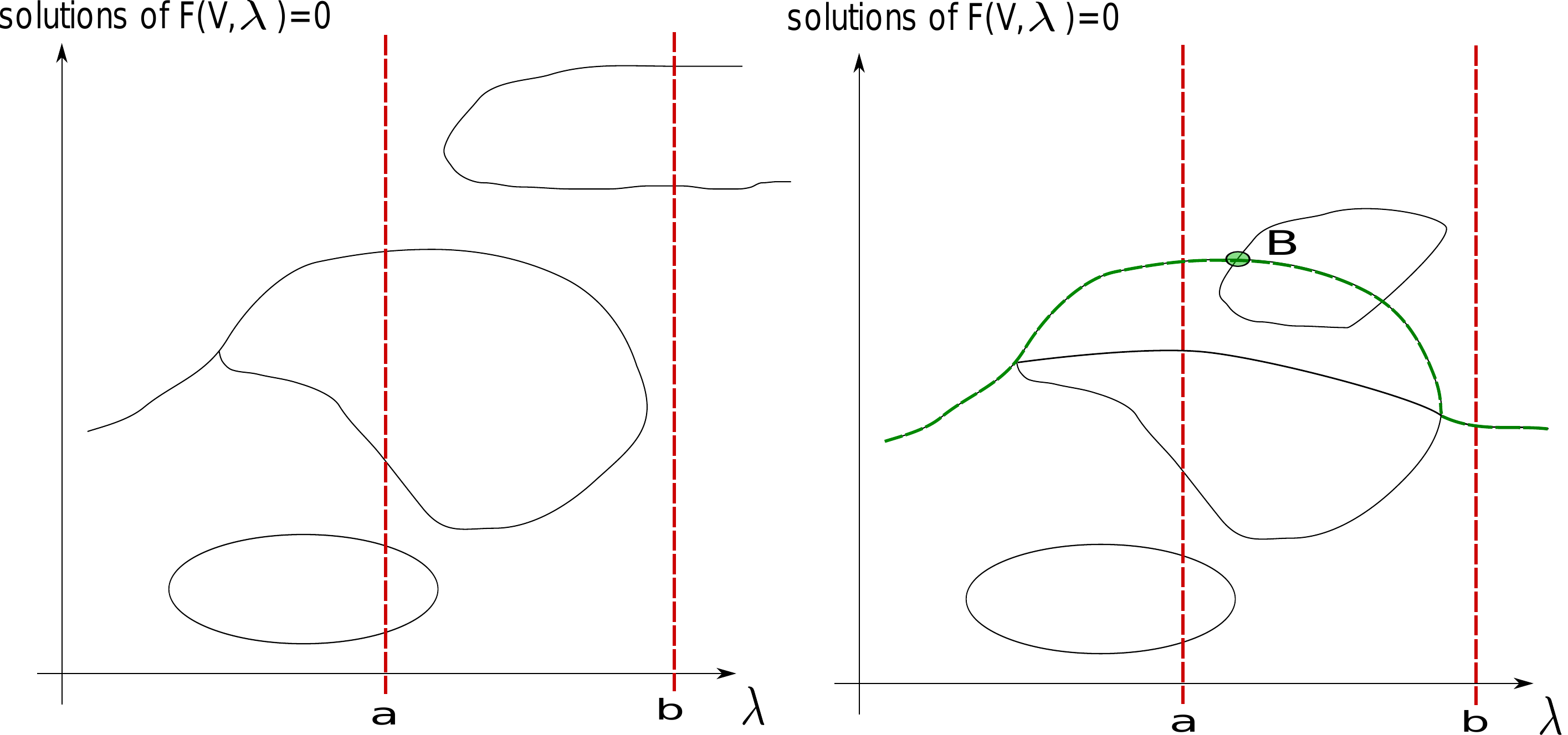}
       \caption{In the lefthand part of the figure, there is no connected curve of solutions in $[a,b]$: this is forbidden by proposition \ref{prop:bumps}.5 which states that we must be in the situation shown in the righthand part of the figure, where the green curve connects $\mathcal{B}_{a} \times \{a\}$ and  $\mathcal{B}_{b} \times \{b\}$, see text.}
\label{fig:CompoConnexe}
\end{figure}

 In fact, except at points where the Jacobian of $F$ is non-invertible (such points like B in figure \ref{fig:CompoConnexe} are potential bifurcation points), the implicit functions theorem tells that $\lambda\to(\vf,\lambda)$ is differentiable. Hence in effect proposition \ref{prop:bumps}.5 imposes strong constraints on the set $\mathcal B_\lambda$ as shown in figure \ref{fig:CompoConnexe}. The horizontal axis represents the parameter $\lambda$, the vertical axis the space $\mF$ where the solutions of (\ref{eq:bump}) live. The curves represent possible solutions as functions of $\lambda$. The configuration in the lefthand part of the figure are forbidden by proposition \ref{prop:bumps}.5 while those on the righthand side are allowed, the green curve being an example of a continuous curve $s\to(\mV^f_{\lambda(s)},\lambda(s))$ from $[0,1]$ to $\mF\times [a,b]$.

\

Proposition \ref{prop:bumps}.5 gives a very interesting general (non-local) property of the set of solutions. But it is non-constructive, for example it does not tell us which branch to chose at point B in figure \ref{fig:CompoConnexe} and we need to compute all the branches to know the path to $\lambda=b$. Hopefully such branching points as B are very rare and one only sees turning points rather than branching points (such as the one at $\lambda_1$ in figure \ref{fig:SigBorne}, left): any perturbation will indeed destroy such branching points (see figure \ref{fig:SigBorne} right). Hence, if one continues the trivial solution (obtained for $\lambda=0$), one will typically find a curve like the green one in figure \ref{fig:SigBorne}.Right. This may lead to the wrong conclusion  that for $\lambda$ big enough there is only one stationary solution instead of three. The problem is to find a way to compute, if it exists, the second, red, curve which is not connected to the green curve, hence not attainable by $\lambda$-continuation.
 
 An idea, directly suggested by the above picture is to restore the branching points by perturbation. Among all possible perturbations, we choose one of the simplest, we vary the contrast $\varepsilon\mI_{\rm ext}$ by varying $\varepsilon$ from 0 to 1, as suggested by experimentalists who usually provide neural responses as functions of the contrast.
 
The conclusion is that if we want one solution for $\lambda\neq0$, we use a $\lambda$-continuation of the trivial solution,  but if we want more than one (or the maximum number of) solutions for $\lambda\neq0$, then we perform a $(\lambda,\epsilon)$-continuation of the trivial solution.

\begin{remark}
An interesting question is to predict how close to $\lambda_L$ can be the smallest value of $\lambda$ where a ``turning point'' occurs. 
\end{remark}

However, when performing this $(\lambda,\varepsilon)$-continuation of the trivial solution, we may end-up with a lot of data. Hence, we need to have an idea of of the structure of the set $(\mV^f_{\lambda,\epsilon},\lambda,\epsilon)$. The idea is, again, to think of this set as a deformation of an easier to compute set of stationary solutions. This is  done in the next section.

\begin{remark}
 It is highly possible that this $(\lambda,\varepsilon)$-continuation scheme still misses some solutions. One possibility is  to introduce a third parameter, but it may become quickly numerically untractable.
\end{remark}

\section{Exploring the set of persistent states}\label{section:bifurcations}

We exploit the fact that in the neural field equation \eqref{eq:bump} the ratio between the external current $\mI_{ext}$ and $\mJ$ is not fixed a priori. Hence when studying the NFE, one would rather look at 
\begin{equation}
\label{eq:net}
-\mV+\mJ \cdot \mS(\lambda\mV)+\varepsilon\mI_{\rm ext}=0 
\end{equation}
where $\varepsilon \geq 0$. The persistent states now depend upon the pair $(\lambda,\varepsilon)$ \textit{i.e.} $\mV^f_{\lambda,\varepsilon}$. The idea is to infer the persistent states of \eqref{eq:net} from those, $\mV_{\lambda,0}$ of 
\begin{equation}
\label{eq:net2}
- \mV+\mJ \cdot \mS(\lambda\mV)=0
\end{equation}
We further simplify the problem by considering
\begin{equation}
\label{eq:net3}
-\mV+\mJ \cdot \mS_0(\lambda\mV)+\mu \mJ\cdot\mS(0)=0,
\end{equation}
where $\mS_0$ is defined in proposition \ref{prop:bumps}.
We recover equation \eqref{eq:net2} when $\mu=1$. The advantage is that when $\mu=0$, we can say a great deal about the persistent states. 
\subsection{A simpler case}\label{section:simpler}
This simpler case reduces to the study of the previous equation when $\mu=0$:
\begin{equation}
\label{eq:Simpler}
- \mV+\mJ \cdot \mS_0(\lambda\mV)=0
\end{equation}
In other words, we infer the persistent states of \eqref{eq:net} from those, $\mV_{\lambda,\varepsilon=0,\mu=0}$ of equation \eqref{eq:Simpler}.
This case has been studied a lot by several authors \cite{ermentrout:98,coombes:05,bressloff:05} when a constant current $\mI_{ext}$ is applied, which amounts to changing the threshold $\boldsymbol{\theta}$ in $\mS$.

$\mV=0$ is a trivial persistent state for \eqref{eq:Simpler}. Recall that a necessary condition for the equation $F(\mV,\lambda)=0$ to bifurcate at the solution $(\mV_{\lambda}^f,\lambda)$ is that $D_VF(\mV^f_\lambda,\lambda)$ is non-invertible.

For technical reasons it is interesting to slightly modify the nonlinear term $\mR$ in \eqref{eq:ChauchyPb} by subtracting from it its linear part as follows.
We rewrite (\ref{eq:Simpler}) as
\begin{equation}
\label{eq:ChauchyPbBif}
\mL_\lambda \cdot \mV+\mR(\lambda,\mV)=0
\end{equation}
with $\mR(\lambda,\mV)=O(\|\mV\|_\mF^2)$ and 
\begin{equation}
\mL_\lambda=-{\rm Id}+\lambda \mJ  D\mS_0(0)=-{\rm Id}+\frac{\lambda}{4} \mJ
\label{eq:Llambda}
\end{equation}
The operator $\mL_\lambda$ satisfies the following properties.
\begin{proposition}\label{prop:fredholm}\
\begin{itemize}
\item The operator $\mJ$ is a compact linear operator on $\mF$.
\item $\mL_\lambda=-{\rm Id}+\frac{\lambda}{4} \mJ$ is a Fredholm operator of index 0 (useful for static bifurcation, see \cite{haragus-iooss:09}).
\item $\mL_\lambda$ is a sectorial operator (useful for dynamical bifurcation, see \cite{haragus-iooss:09}).
\end{itemize}
\end{proposition}
\begin{proof}
The first property was proved in lemma \ref{lemma:R}.
 
Because $\mJ$ is a compact operator, the kernel of $-{\rm Id}+\frac{\lambda}{4}\mJ$ is of finite dimension equal to the codimension of its image, hence its index is 0.

$\mL_\lambda$ is a sectorial operator because it is a bounded linear operator on a Banach space, see e.g. \cite{henry:81}.
\end{proof}

\

We can therefore state that the values of the parameter $\lambda$ that determine the \textit{possible} bifurcation points are:
\[
 \lambda_n=\frac{4}{\sigma_n},\,n=1,\,2,\cdots
\]
where $\sigma_n$ is an eigenvalue of the compact operator $ \mJ$. We assume in the sequel that $\lambda_1 \leq \lambda_2 \leq \cdots$.

Additional properties have to be met in order to have a bifurcation point (see \cite{ma-wang:05}) but they are satisfied if every eigenvalue of $ \mJ $ is simple. We will assume that the first $k\geq 1$ ($k$ is arbitrary) eigenvalues are simple because this class of operators is dense in the set of compact operators set (see \cite{ma-wang:05}). We denote by $e_n$ (respectively by $e_n^*$) the eigenvector of $ \mJ $ (respectively of the ajoint operator $\mJ^*$) for the eigenvalue $\sigma_n$. 

A simple adaptation of lemma \ref{lemma:R} shows that $\mR$ is $C^q$, for all integer $q$ and we can consider (if it exists the minimal integer $q\geq 2$ such that 
\[
0\neq\chi_q^{(n)}\equiv \left\langle  G_q(e_n,\lambda_n),e_n^*\right\rangle_\mF.
\]

According to lemma \ref{lemma:R} 
\[
 G_q(e_n,\lambda_n)=\frac{1}{q!}D^q \mR(\lambda_n,\mathbf 0)[e_n \cdots e_n]=\frac{\lambda_n^q}{q!}S_0^{(q)}(0) \mJ \cdot (e_n\cdots e_n)
\]

By parity, $S^{(q)}_0(0)\overset{\rm def}{=} s_q \neq 0$ if and only if $q$ is odd. Hence, the parity of $q$, which is odd, tells that we have a pichfork bifurcation at $\lambda_n$. In particular, the bifurcated branch is given by

\[
 \mV_\lambda^f=x(\lambda) e_n + o(x)
\]

Thus the bifurcation portrait is a set of branches $\mathcal C_n$ emanating at points (at least for $n\leq k$) $(0,\lambda_n)$. Thanks to proposition.\ref{prop:fredholm}, we can apply the Lyapunov-Schmidt procedure (see \cite{haragus-iooss:09}) to get the following reduced equation on $x\in\mathbb R$ :
\[
 0=(-1+\lambda\sigma_n/4) x+g(x,\lambda)=(-1+\lambda\sigma_n/4) 
x+\chi^{(n)}_q x^q+o(x^q)
\]

Depending on the sign of $\chi_q^{(n)}$, the pichfork branch is oriented toward $\lambda<\lambda_n$ (resp. $\lambda>\lambda_n$) if $\chi^{(n)}_q>0$ (resp. $\chi^{(n)}_q<0$).
Let us note $e_n^q$ the vector $\underset{q \rm times}{\underbrace{e_n \cdots e_n}}$. We have
\[
 \chi_q^{(n)}=\frac{\lambda_n^q s_q }{q!} \left\langle \mJ \cdot e_n^q, e_n^* \right\rangle_\mF = \frac{\lambda_n^q s_q }{q!} \left\langle  e_n^q, \mJ^* \cdot e_n^* \right\rangle_\mF=\frac{\lambda_n^{q-1} s_q }{q!s_1} \left\langle e_n^q, e_n^* \right\rangle_\mF
\]
From a practical standpoint, the last inner product is difficult to compute since it requires the computation of the eigenvectors of the adjoint operator $\mJ^*$ of $\mJ$ for the inner-product of the Hilbert space $\mF=\mW^{m,2}$ which are in general different from those of the adjoint operator $\mJ^*_2$ of $\mJ$ for the inner-product of
the Hilbert space $\mW^{0,2} \equiv \mL^2$. This computation is greatly simplified by the following proposition.
\begin{proposition}\label{prop:inner}
 Let $e^*_{L^2}$ be an eigenvector of the adjoint operator $\mJ^*_{L^2}$ of $\mJ$ for the inner-product of
the Hilbert space $\mL^2$, associated to the eigenvalue $\lambda$ (assumed to be simple). Let $e^*_\mF$ be the corresponding eigenvector of the adjoint operator $\mJ^*_\mF$ of $\mJ$ for the inner-product of the Hilbert space $\mF=\mW^{m,2}$. Then it is possible to choose $e^*_\mF$ such that the following holds
\[
 \left\langle \mU, e^*_\mF \right\rangle_\mF=\left\langle \mU, e_{L^2}^* \right\rangle_2 \quad \forall \mU \in \mF
\]
\end{proposition}
\begin{proof}
see appendix.\ref{appendix:ajoint}.
\end{proof}

\

In effect, in order to obtain the sign of $\chi^{(n)}_q$, we only have to compute the eigenvectors of $\mJ_2^*$ and to compute inner products in $\mL^2$.

\

Thus, we have found local branches of stationary solutions and continue them globally in order to obtain the global branches named $\mathcal C_n$. An interesting question, yet unsolved, is to know whether the branches $\mathcal C_n$ are connected. Some results exist in that direction in the line of those of Rabinowitz (see \cite{rabinowitz:71,ma-wang:05}) but dont give  much insights in the general case ($d>1,p>1...$). However, they can be used to derive some properies of $\mathcal C_1$.

\begin{proposition}[Turning point property]
\label{prop:saddle}
If $\chi^{(n)}_q>0$, then 
 $\exists\lambda_T<\lambda_1$ such that $\forall\lambda\in(\lambda_T,\lambda_1)$, (\ref{eq:Simpler}) has at least 3 solutions and $\lambda_T=\min\left\lbrace \lambda/(\mV^f,\lambda)\in\mathcal C_1\right\rbrace $.
\end{proposition}
\begin{proof}
Let $\mathcal C_1$ be the connected component in $\bar B$ where $B=\left\lbrace(\mV,\lambda)|\mV\neq \textbf{0},\ (\mV,\lambda) \text{ satisfies } (\ref{eq:Simpler}) \right\rbrace$   to which $(\textbf{0},\lambda_1)$ belongs. Then (see \cite{ma-wang:05})
$\mathcal C_1$ is unbounded in $\mF\times\mathbb R_+$ or contains a point $(\textbf{0},\lambda_n),\ n>1$. In either case, $\mathcal C_1$ exists until $\lambda=\lambda_2$.

$\tilde{\mathcal C}=\bar{\mathcal C_1}\cap(\mF\times[\lambda^*,\ \lambda_2])$ is closed and bounded (because every solution $\mV^f_\lambda$ is bounded in $\mF$). Let us show that $\tilde{\mathcal C}$ is compact in $\mF\times [0,\ \lambda_2]$. 
Consider a sequence $(\mV_n,s_n)$ in $\tilde{\mathcal C}$. As $s_n$ is bounded, we can assume it is convergent. We also have  $\mV_n=\mJ \cdot\mS(s_n\mV_n)+\mI_{\rm ext}$. As $(\mV,\lambda)\to\mJ \cdot \mS(\lambda\mV)+\mI$ is a compact operator, there exists a subsequence $(\mV_{\phi(n)},s_{\phi(n)})$ such that $\mJ \cdot \mS(s_{\phi(n)}\mV_{\phi(n)})+\mI$ is convergent, hence $\mV_{\phi(n)}$ is converging. We have proved that $\tilde{\mathcal C}$ is compact. 
Hence $\Pi_{\mathbb R^+}(\tilde{\mathcal C})$ is a compact subset of $\mathbb R^+$. Then $\inf\Pi_{\mathbb R^+}(\tilde{\mathcal C})$ is a $min$ written $\lambda_T\in\Pi_{\mathbb R^+}(\tilde{\mathcal C})$. As it is an $\inf$, there exists a sequence $s_n$ associated to a $\mV_n$ in $\tilde{\mathcal C}$ such that $s_n\to\lambda_T$. But as $\tilde{\mathcal C}$ is compact, we can assume that $\mV_n\to\mV_T$. Then $(\mV_T,\lambda_T)\in\mathcal C$ is called a \emph{turning point}\footnote{Not in the sense of \cite{kuznetsov:98}, here it denotes a local extremum in the parameter along the curve of solutions.}.

So we have proved that $\lambda_T\leq \lambda_1$. But in the case $q$ odd with $\chi^{(n)}_q>0$, $\mathcal C_1$ exists for $\lambda<\lambda_1$ and $\lambda_T<\lambda_1$.

Now, $\forall\lambda\in(\lambda_T,\lambda_1),\ \exists (\vf,\lambda)\in\mathcal C_1\setminus\left\lbrace (\textbf{0},\lambda)\right\rbrace$ because $\mathcal C_1$ intersects $\left\lbrace \textbf{0}\right\rbrace\times\mathbb R_+$ only at bifurcation points $(\textbf{0},\lambda_n)$ located 'after' $\lambda_1$. Then because of  proposition \ref{prop:bumps} part 4, there are at least three solutions.
\end{proof}

\begin{remark}
 If the sigmoidal function had satisfied $S^{(2)}(0)\neq 0$, we would have seen transcritical bifurcations and the previous proposition still holds in that case.
\end{remark}

\begin{remark}
 If we were able to prove that the branches do not intersect, the previous proposition would apply to all branches $\mathcal C_n$ satisfying the required conditions.
\end{remark}

Even if we do not deal with the dynamics, we can say a little using  \cite{haragus-iooss:09,ma-wang:05,kuznetsov:98}. At every bifurcation point, the center manifold at $\mV=0$ (see \cite{haragus-iooss:09}) is one-dimensional while the dimension of the unstable manifold increases as $\lambda$ crosses values corresponding to transcritical points. Hence for large $\lambda$s, we can have a large unstable manifold. Note that every value of $\lambda$ is biologically plausible because the locally (around $\mV=0$) exponentially divergent dynamic is bounded (see proposition \ref{prop:bounded}), which can make the ``global'' dynamics \textbf{very} intricate. In effect, when $\lambda$ grows to infinity, the sigmoid tends to a Heaviside function which acts as a threshold.

\begin{remark}
 Note that it is easy to characterize Breathers (see \cite{nishiura-mimura:89,coombes-owen:05,folias-bressloff:05}) in this framework. It is sufficient to choose a connectivity function $\mJ$ with complex eigenvalues. For example with a model of two populations ($p=2$), one excitatory, one inhibitory, such that:
 \[\mJ(\mr,\mr')=
 \left[ 
\begin{array}{cc}
  G_1(\mr,\mr') &  - G_2(\mr,\mr')\\
   G_2(\mr,\mr') &  G_1(\mr,\mr')
\end{array}
\right] \quad G_i\text{ symetric} \geq 0,\ [G_1,G_2]=0 
\]
Then $\sigma(\mJ)=\left\lbrace   g_{1,n}\pm i g_{2,n},g_{i,n}\in\sigma(G_i)\right\rbrace $ and a Hopf bifurcation occurs at $\lambda_{  H,n}=\frac{1}{2S'(0)g_{1,n}}$ if the eigenvalue is simple. Because $\mR$ is $C^3$, we can also compute its first Lyapunov coefficient (see \cite{haragus-iooss:09},\cite{kuznetsov:98}{p.107}) and find 
\[
  l_1=-\frac{\lambda_{H,n}^3s_2^2}{4\omega_n^2s_1}\left( 1-\sqrt{g_{1,n}^2+g_{2,n}^2}\frac{s_3s_1}{s_2^2}\right) 
\]
with   $\omega_n=\frac{\sqrt{g_{1,n}^2+g_{2,n}^2}}{2g_{1,n}},\ 
	s_2=S^{(2)}(0),\
  	s_3=S^{(3)}(0)$
\end{remark}

\subsubsection{The case $d=1,p=1$}\label{section:rabinowitz}

We can use the results of Rabinowitz (see \cite{rabinowitz:71}) in the case $d=p=1$ and when the connectivity operator $\mJ$ is a symmetric convolution. In that case the zeros of the eigenstates of $\mJ$ are simple on $\Omega\subset\R$ (because they are cos and sin functions). Moreover the Taylor formula allows to write $\mS_0(x)=xh(x)$ with $h>0$. Then from \cite{rabinowitz:71}, it follows that the bifurcated branches\footnote{Here they are connected components in $C^1(\Omega,\R)\times\R_+$ rather than in $\mF\times\R_+$.} $\mathcal C_n$ (pichfork or transcitical)  meet $(\lambda_n,0)$ and $\infty$ and are caracterized by the number of zeros of their elements ; namely $\forall\ (\mV^f,\lambda)\in\mathcal C_n$, $\mV^f$ has exactly $n-1$ zeros in $\Omega$. As a consequence, the branches $\mathcal C_n$ do not intersect.

If we were able to prove that all the stationary solutions are connected to the zero solution, we would have completely caracterized $\mathcal B$. Nevertheless, it is still a lot of information.

\begin{remark}
 The study of the simpler case is very important for the numerical computations. Indeed, it allows analytical predictions. When one performs the $(\lambda,\mu,\epsilon)$-continuation, one should look at the section of solutions $(\lambda,\mu=0,\epsilon=0)$ and compare to the predictions of the simpler case to see if the numerical did not miss some solutions (that may happen when the system has some symmetries).

\end{remark}

\subsection{Returning to the original equation}

The overall picture that emerges from the previous section is interesting despite the fact that some of its features are hard to justify from the biological viewpoint: no external input, rate function $\mS_0$ possibly negative. The reason is that it gives clues about the set of persistent states when $\varepsilon \neq 0$ and provides a way to compute them numerically. Starting with the trivial solution $\mV^f_0$ of (\ref{eq:net2}) when the slope parameter $\lambda$ is null, we can perform a numerical continuation with respect to the two parameters $\lambda,\ \varepsilon$. 

When $\varepsilon \neq 0$, the only bifurcations that are possibly unaltered are the turning points. The transcritical/pichfork bifurcations will be ``opened'' as described below. We are still able to predict the stability near the points $(\mathbf{0},\lambda_n)$.
Let us note  $\bar I=\left\langle \mu \mJ\cdot\mS(0)+\varepsilon\mI_{ext},e_n^*\right\rangle_\mF $. Then, the Lyapunov-Schmidt  method \cite{haragus-iooss:09,ma-wang:05}  leads to 
\begin{equation}
 \label{eq:LSbarI}
0=(-1+\frac{\lambda\sigma_n}{4})x+\chi^{(n)}_qx^q+\bar I+o(x^q) 
\end{equation}

\begin{figure}[htbp]
 \centerline{
          \includegraphics[width=0.5\textwidth]{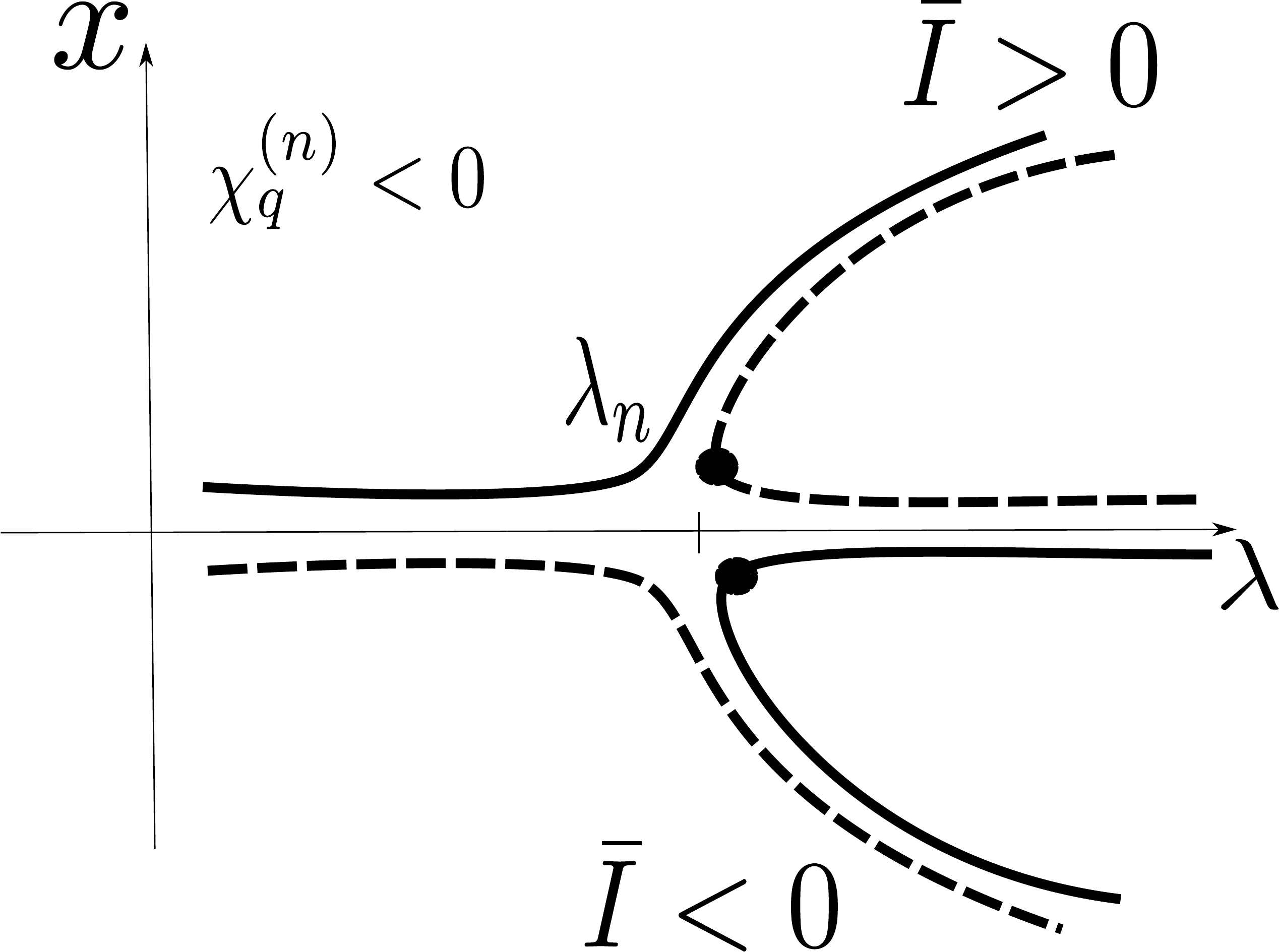}
	  \includegraphics[width=0.5\textwidth]{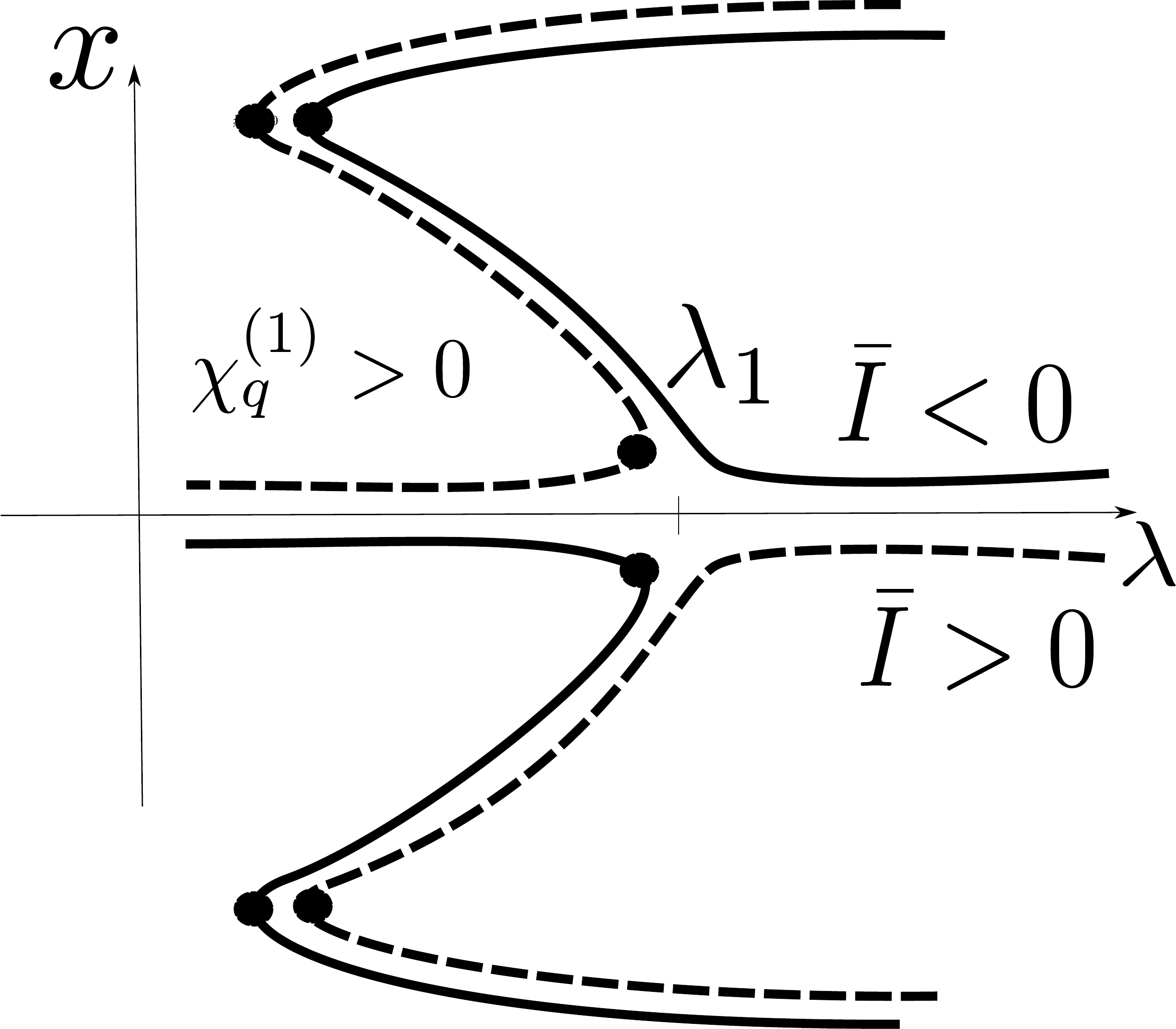}
}
          \caption{Opening of the pichfork bifurcation at the first (right) and $n$th ($n>1$, left) eigenvalue.
Black dots indicate saddle nodes. Note that there are three such points, see proposition \ref{prop:saddle}, for the first eigenvalue. Note that the branch may have a more intricate shape than the one shown.}
\label{fig:opening}
\end{figure}

Notice that when $q=2$, we find the same equilibria as those of the Bogdanov-Takens normal form and for $q=3$, we obtain a cusp.
Solving the  polynomial equation (\ref{eq:LSbarI}) allows us to describe different opening scenarios depending on the sign of $\bar I$. The cases $\bar I>0$ are shown on figure \ref{fig:opening}. The case of $\lambda_1$ is a little bit special according to proposition \ref{prop:saddle} and is shown in the righthand part of figure \ref{fig:opening}.

\

\begin{remark} 
When the first eigenvalue generates a 'well oriented' pichfork branch, proposition \ref{prop:saddle} says that a turning point must occurs on this branch: there are two turning points on this branch. Still we have a local description and it would be interesting to have more global results for example concerning the behaviour of these turning points when $(\lambda, \varepsilon)$ varies. 
\end{remark}

\section{Reduction  to a finite dimensional analysis}\label{section:PG}

Neural field models are one of the possible generalizations of standard neural networks considered as discrete sets of connected neurons. They can be characterized by two limit processes. First, we let the total number of neurons grow to infinity so that each node of the network represents an ideally infinite number of neurons, in practice a large number of such neurons belonging to different populations, in effect a neural mass. Second, we assume that these neural masses form a continuous neural material and let the connectivity graph of the neural network become continuous. The graph connectivity matrix then turns into a continuous function of the spatial coordinates.  One would think that after passing twice to the limit the resulting system would be infinite-dimensional. However, the dimensionality of the neural field models depends essentially upon the linear operator $\mJ$ representing the connectivity function. If this linear operator has a finite-dimensional range, we show below that the corresponding neural field model is finite-dimensional and is equivalent to a finite set of ordinary differential equations. Moreover, even if this condition is not met we also show that we can always approximate the operator as accurately as desired by a finite-dimensional range operator and reduce the neural field model to a finite set of ordinary differential equations.

\subsection{The Pincherle-Goursat Kernels}
In our numerical studies we use connectivity functions that are such that the corresponding linear operators of $\mF$ have finite rank, i.e. their range is a finite dimensional subspace of $\mF$. This is without loss of generality because of the following theorem (see, e.g. \cite{brezis:83}):

\begin{theorem}
 The subspace $\mcR_f(\mF)$ of finite-dimensional range linear operators of $\mF$  is dense in $\maH$, the set of linear compact operators of $\mF$.
\end{theorem}
\begin{proof}
 This is true because $\mF$ is a Hilbert space \cite[Chapter 6]{brezis:83}
\end{proof}

\noindent
In the area of integral equations, these operators are called Pincherle-Goursat kernels \cite{tricomi:85}, in short PG-kernels. 
They are defined as follows

\begin{definition}[Pincherle-Goursat Kernels]
The connectivity kernel $\mJ(\mr,\mr')$ is a PG-kernel if 
\[
\mJ(\mr,\mr')=\sum\limits_{k=1}^{N}X_k(\mr)\otimes Y_{k}(\mr')
\]
where $X_k,\,Y_k,\,k=1\cdots N$ are two sets of $N$ linearly independent elements of $\mF$, and $X_k(\mr)\otimes Y_{k}(\mr')$ is the rank 1 $p \times p$ matrix $X_k(\mr)Y_{k}^T(\mr')$.
\end{definition}

\noindent
We have 
\[
\mJ \cdot \mU=\suml_{k} X_k \langle Y_k, \mU \rangle_2,
\]
where $\langle \cdot,\cdot \rangle_2$ is a short notation for the inner-product in $\mL^2(\Omega,\R^p)$. Thus $\mJ \cdot \mU\in {\rm Span}(X_1,\cdots,X_{N})$ which we denote by $R(\mJ)$. 

\subsection{Persistent state equation for PG-kernels}
We now cast the problem of the computation of the solutions of equation \eqref{eq:bump} into the PG-kernel framework:
\[
\mV-\mI_{\rm ext}=\mJ \cdot \mS(\lambda \mV)
\]
Since $\mV-\mI_{\rm ext}\in R(\mJ)$, we can write $\mV-\mI_{\rm ext}=\suml_k v_k X_k$, and note $v=(v_k)_{k=1\cdots N}$.
The persistent state equation reads:
\[
v_k=\left\langle Y_k,\,\mS\left(\lambda\left( \suml_k v_k X_k+\mI_{\rm ext} \right)\right) \right\rangle_2 \quad k=1,\cdots, N
\]
This is a set of $N$ nonlinear equations in the $N$ unknowns $v_1,\cdots,v_N$ which can be solved numerically using classical methods. 

\subsection{Reduction to a finite number of ordinary differential equations}
In this section we reduce equation \eqref{eq:NM} (when $\mJ$ is independant of $t$) to a system of ODEs. We write $\mI$ instead of $\mI_{\rm ext}$ for simplicity and note $\mS_\lambda$ the function $\R^p \to \R^p$ defined by $\mS_\lambda(x)=\mS(\lambda x)$ for $x \in \R^p$.

We note $R(\mJ)^\bot$ the orthogonal complement of $R(\mJ)$  in $\mF$ with respect to the inner-product $\left\langle ,\right\rangle_2$:
\[
 \mF=R(\mJ)\oplus R(\mJ)^\bot
\]
We write
\[
 \mV=\mV^\parallel+\mV^\bot,
\]
where $\mV^\parallel$ (resp. $\mV^\bot$) is the orthogonal projection of $\mV$ on $R(\mJ)$ (resp. $R(\mJ)^\bot$). We have a similar decomposition for the external current $\mI$
\[
 \mI=\mI^\parallel+\mI^\bot
\]
We now decompose $R(\mJ)$ as a Cartesian product of $p$ finite dimensional subspaces of $\mG$.
Because 
\[
 \mJ_{ij}(\mr,\mr')=\sum_{k=1}^NX_k^i(\mr)Y_k^j(\mr') \quad i,\,j=1,\cdots,p,
\]
each coordinate $V_i$, $i=1,\cdots,p$ of $\mV$ satisfies 
\[
\dot V_i+\alpha_i V_i=\sum\limits_{k=1}^N \left\langle Y_k,\mS(\lambda\mV)\right\rangle_2 +I_i \quad i=1,\cdots,p
\]
 Let us consider the $p$ finite dimensional subspaces $E_i$, $i=1,\cdots,p$ of $\mG$, where each $E_i$ is generated by the $N$ elements  $X^i_k,\,k=1,\cdots,N$.  We note $E_i^\bot$ the orthogonal complement of $E_i$ in $\mG$. This induces a decomposition of $\mF$ as the direct sum of the cartesian product $\prod_{i=1}^p E_i=R(\mJ)$ and its orthogonal complement $\prod_{i=1}^p E_i^\bot = R(\mJ)^\bot$. We write $V_i=V_i^\parallel+V_i^\bot$ as well as $I_i=I_i^\parallel+I_i^\bot$. We then have 

 \begin{equation}
  \begin{cases}
   \dot{V}_{i}^\parallel+\alpha_i V_{i}^\parallel=\sum\limits_{k=1}^N \left\langle Y_k,\mS(\lambda\mV)\right\rangle_2 +I_{i}^\parallel\\
   \dot{V}_i^\bot+\alpha_i V_i^\bot=I_i^\bot
  \end{cases}i=1,\cdots,p
 \end{equation}
Considering the canonical basis $e_i$, $i=1,\cdots,p$, of $\R^p$, we define 
\[
 \begin{array}{ll}
  	\mV^\parallel=\sum\limits_{i=1}^p V_i^\parallel e_i & \mI^\parallel=\sum\limits_{i=1}^p I_i^\parallel e_i\\
	\mV^\bot=\sum\limits_{i=1}^p V_i^\bot e_i & \mI^\bot=\sum\limits_{i=1}^p I_i^\bot e_i
 \end{array}
\]
We obtain the $2p$-dimensional non-autonomous system of ODEs:
\begin{equation}
  \begin{cases}
   \dot\mV^\parallel+\mL \cdot \mV^\parallel=\mJ\cdot\mS(\lambda\mV) +\mI^\parallel\\
   \dot\mV^\bot+\mL \cdot \mV^\bot=\mI^\bot \\
  \end{cases}
 \end{equation}

\begin{remark}
 If $\mI^\bot$ is stationary then $\mV^\bot$ converges to $\mL^{-1}\mI^\bot$.
\end{remark}

\section{One population of orientation tuned neurons: the ring model}\label{section:ring}

As an application of the previous results, we study the  ring model of orientation tuning introduced by Hansel and Sompolinski (see \cite{hansel-sompolinsky:97,shriki-hansel-etal:03,ermentrout:98,dayan-abbott:01,bressloff-bressloff-etal:00,bressloff-cowan-etal:01}), after the work of Ben-Yishai (see \cite{ben-yishai-bar-or-etal:95}), as a model of a hypercolumn in primary visual cortex. It can be written as:
\[
   \tau\dot A(x,t)= -A(x,t) + S\left(\lambda \left(\int\limits_{-\pi/2}^{\pi/2}  J(x-y)  A(y,t))dy/\pi+I(x)-\theta\right)\right)
\]

Some authors, \cite{bressloff-bressloff-etal:00,bressloff-cowan-etal:01}, chose $J$ to be a difference of Gaussians. On the other hand, Ben-Yishai, in \cite{ben-yishai-bar-or-etal:95}, started with a network of excitatory/inhibitory spiking neurons and derived a meanfield approximation of this network yielding the activity response described by the following equations:
\[
 \begin{cases}
   \tau\dot A(x,t)= -A(x,t) + S\left(\lambda \left(\int\limits_{-\pi/2}^{\pi/2}  \left[J_0+J_1 \cos(\alpha(y-x))\right]  A(y,t))dy/\pi+\varepsilon I(x)-\theta\right)\right)\\
   I(x) =  1-\beta+\beta \cos(\alpha(x-x_0))\\
  \end{cases}
\]
$\alpha=2$, $0\leq \beta \leq 1$ and the threshold $\theta=1$  in the above cited papers.
Being an activity model and not a voltage model in the terminology of \cite{ermentrout:98,faugeras-grimbert-etal:08} it is not directly amenable to our analysis. We can either extend this analysis to activity models as shown in appendix \ref{appendix:activity} or do the following. We rewrite the previous equation as
\[
 \tau\dot A=-A+S(\lambda(J \cdot A+\varepsilon I-\theta)),
\]
and perform the change of variable $V=J \cdot A+ \varepsilon I-\theta$. This leads to the following equations
\begin{equation}
  \begin{cases}
  \label{eq:RM3d}
   \tau\dot V(x,t)= -V(x,t) + \int\limits_{-\pi/2}^{\pi/2}  \left[J_0+J_1 \cos(\alpha(y-x))\right]  S(\lambda V(y,t))dy/\pi+\varepsilon I(x)+\theta\\
   I(x) = 1-\beta+\beta \cos(\alpha(x-x_0))
  \end{cases}
 \end{equation}
We are now in the case of the model studied in this paper with $p=1$, $d=1$ and $\Omega=(-\pi/2,\pi/2)$. Note that since the $x$-coordinate represents an angle, this is not a neural field model per se but rather a neural mass model.

The nonlinearity is  often chosen to be a Heaviside function, or, as in \cite{ben-yishai-bar-or-etal:95}, a piecewise linear approximation of the sigmoid \footnote{This does not allow for computation of bifurcation branches but allows the detection of branching points.}, or, as in \cite{ermentrout:98,dayan-abbott:01,bressloff-bressloff-etal:00}, a true sigmoidal function. $J_1$ can take any values and $I$ is an external current coming from the LGN. $J_0$ is most of the time negative (see \cite{ben-yishai-bar-or-etal:95,dayan-abbott:01,bressloff-bressloff-etal:00,bressloff-cowan-etal:01}) but can be positive as well(see \cite{bressloff-cowan-etal:01}): the $J_i$s can be thought of as the first Fourier coefficients of $J$, $J_0$, being its mean value, it can be positive even if the surround is inhibitory. We can, up to a multiplication of the previous equation, make the assumption 
\[
J_0=\varepsilon_0\in\left\lbrace-1,1 \right\rbrace 
\]
 
 For example, in \cite{dayan-abbott:01}, we find $J_0=-7.3,\ J_1=11,\ \beta=0.1,\ \theta=0$ which are taken from \cite{ben-yishai-bar-or-etal:95} except for $\theta=1$. The slope is assumed to be $\lambda=1$. Using the previous scaling, it becomes $J_0=-1,\ J_1=1.5,\ \lambda=7.3/s_1=29.2$ and $\theta\to\theta/7.3$ which gives $\theta\approx 0.1$ in the case of \cite{ben-yishai-bar-or-etal:95} and $\theta=0$ in \cite{dayan-abbott:01}.
 
\
 
The goal of this section is not to derive the whole bifurcation diagram of the Ring Model but rather to show how the stationary solutions are organized and to give clues about the dynamics in a given range of parameters. This study is helpful because many large scale model of V1 use the Ring Model for the hypercolumns. We will see that, depending on the stiffness of the nonlinearity, there may  exist many stationary solutions, which are all acceptable responses of the network for a given input of the LGN. Thus these local orientation detectors may behave less trivially than they were initially made for. These cortical states can make the local dynamics sophisticated when the slope parameter $\lambda$ is big enough.

\subsection{Mapping the ring model to the PG-kernel formalism}
 Expanding the cosine in the previous equation, and denoting by $\cos_\alpha$ (respectively $\sin_\alpha$) the function $x \to \cos(\alpha x)$ (respectively $x \to \sin(\alpha x)$), we find that, depending on the sign of $J_i$, ($\varepsilon_i=sign(J_i),i=0,1$): 
\[ 
 J = \varepsilon_0 1\otimes 1+\varepsilon_1\sqrt{|J_1|}\cos_\alpha\otimes\sqrt{|J_1|} \cos_\alpha+ \varepsilon_1\sqrt{|J_1|} \sin_\alpha\otimes \sqrt{|J_1|}\sin_\alpha =\sum\limits_{i=0}^2\varepsilon_i X_i\otimes X_i,\,\varepsilon_2=\varepsilon_1.
\]
This formulation has the advantage of preserving the symmetries of $\mW$. With the notations of the previous section, we have $I^\bot=0$, and 
\[
 	\tau\dot V^\parallel=-V^\parallel + W \cdot S(\lambda(V^\parallel+V^\bot_0 e^{-t/\tau}))+  \varepsilon I^\parallel+\theta
\]
where
\begin{equation}\label{eq:Vparallel}
	V^\parallel(x,t)=v_1(t)+v_2(t)\sqrt{|J_1|} \cos_\alpha x+v_3(t)\sqrt{|J_1|} \sin_\alpha x,
\end{equation}

\textit{i.e.} the model is three-dimensional. Note that the previous equation is a rewriting of (\ref{eq:RM3d}) \textbf{without any approximation}.

Similarly we have
\begin{equation}\label{eq:Iparallel}
I=I^\parallel= 1-\beta+ \frac{\beta \cos_\alpha x_0}{\sqrt{|J_1|}} \stackrel{X_1}{\overbrace{\sqrt{|J_1|} \cos_\alpha x}} + \frac{\beta \sin_\alpha x_0}{\sqrt{|J_1|}} \stackrel{X_2}{\overbrace{\sqrt{|J_1|}\sin_\alpha x}}
\end{equation}
As $V^\bot(t)\to 0$, we restrict the study to the case $V^\bot=0$ even if we lose some of the 'real' dynamics by doing so. This is motivated by the fact that the dynamics is made of heteroclinic orbits (as we will see in a moment) between persistent states belonging to the vectpr space $V^\bot=0$. Hence, using this simplification, we are led to study the following 3D system:
\begin{equation}
 \begin{cases}
 \label{eq:RMPG1}
  \tau\dot v_1 = - v_1  +\varepsilon_0 \left\langle S(\lambda V),\text{1}\right\rangle+\varepsilon I_1^\parallel+\theta\\
  \tau\dot v_2 = - v_2  +\varepsilon_1\sqrt{|J_1|}\left\langle S(\lambda V),\cos_\alpha\right\rangle+\varepsilon I_2^\parallel\\
  \tau\dot v_3 = - v_3  +\varepsilon_1\sqrt{|J_1|}\left\langle S(\lambda V),\sin_\alpha\right\rangle+\varepsilon I_3^\parallel
  \end{cases} 
\end{equation}
where $\left\langle f,g\right\rangle=\int\limits_{-\pi/2}^{\pi/2}f(x)g(x)\frac{dx}{\pi}$, $V$ is given by equation \eqref{eq:Vparallel} and $I_i^\parallel$, $i=1,2,3$ is given by equation \eqref{eq:Iparallel}. Note that the basis $(X_0,X_1,X_2)$ is not orthogonal for this inner product.

This system enjoys the  symmetries described by the following lemma.
\begin{lemma}\label{lemma:sym}
When $I_3^\parallel=0$, if $\mv=(v_1\ v_2\ v_3)$ is a solution, then so is $(v_1\ v_2\ -v_3)$.
The plane $v_3=0$ is invariant by the dynamics.
\end{lemma}
\begin{proof}
 This is a consequence of the fact that $\sin_\alpha$ is an odd function while $\cos_\alpha$ is an even function.
\end{proof}

%

\

It is easy to see that 
\[
E(\mv)=-\frac{\|\mv\|^2}{2\tau}+\frac{1}{\tau\lambda}\int_{-\pi/2}^{\pi/2}\bar S_0(\lambda v_1\varepsilon_0+\lambda v_2\varepsilon_1\sqrt{|J_1|} \cos_\alpha x+\lambda v_3\varepsilon_1\sqrt{|J_1|} \sin_\alpha x) \frac{dx}{\pi}+\left\langle\varepsilon I^\parallel+\theta,\mv \right\rangle,
\] 
where $\bar S_0$ is a primitive of $S_0$, is an energy function for the dynamics, i.e. $\dot \mv=\nabla E(\mv)$. Consequently, even for $I$ non spatially homogenous, there are no non-constant periodic trajectories nor homoclinic orbits\footnote{This follows from the consideration of the time derivative of the energy $E$.}. Moreover, all bounded trajectories are stationary solutions or trajectories converging to stationary solutions. Having proven that all trajectories are bounded for the neural field equations in proposition \ref{prop:bounded}, we have characterized the dynamics. It remains to compute the stationary solutions and their attraction basins.

\begin{remark}
 We can generalize these facts to PG-kernels of the type $\mJ=\suml_{i=0}^{N-1}\varepsilon_i X_i\otimes X_i$ by choosing $E(\mv)=-\frac{\left\langle \mv,\mL \mv \right\rangle}{2}+\frac{1}{\lambda}\suml_{k=1}^p\int_\Omega\bar\mS(\suml_i^N \lambda \varepsilon_iv_i X_i^k(\mr))\,d\mr+\left\langle\mI,\mv \right\rangle$ 
\end{remark} 
\subsection{Finding the persistent states}
In order to characterize the set $\mathcal{B}$ of stationary solutions, we apply the scheme of section \ref{section:simpler}. Hence we study the following equation
\[
 V=J\cdot S_0(\lambda V)+\varepsilon I^\parallel+\mu(\theta+J\cdot S(0))
\]
 The nonlinearity is the \textit{odd} function :
 \[
  S_0(x)=\frac{1}{1+e^{-x}}-\frac{1}{2}
 \]
Note that $J\cdot S(0)=\frac{1}{2}(\varepsilon_0+J_1\frac{2sin_\alpha(\pi/2)}{\alpha\pi}cos_\alpha)=\frac{1}{2}\varepsilon_0+\varepsilon_1\sqrt{| J_1|}\frac{sin_\alpha(\pi/2)}{\alpha\pi} X_1$. This gives : 

\begin{equation}
 \begin{cases}
 \label{eq:RMPGBumps}
  \tau\dot v_1 = - v_1  +\varepsilon_0 \left\langle S_0(\lambda V),\text{1}\right\rangle+\varepsilon I_1^\parallel+\mu(\theta+\frac{\varepsilon_0}{2})\\
  \tau\dot v_2 = - v_2  +\varepsilon_1\sqrt{|J_1|}\left\langle S_0(\lambda V),\cos_\alpha\right\rangle+\varepsilon I_2^\parallel+\mu\varepsilon_1\sqrt{ |J_1|}\frac{sin_\alpha(\pi/2)}{\alpha\pi}\\
  \tau\dot v_3 = - v_3  +\varepsilon_1\sqrt{|J_1|}\left\langle S_0(\lambda V),\sin_\alpha\right\rangle+\varepsilon I_3^\parallel
  \end{cases} 
\end{equation}
\subsubsection{The simpler case $\mu=\varepsilon=0$}
This corresponds to finding the persistent states when $\mu=\varepsilon=0$ ensuring that $\mv=0$ is a solution.
For the sake of  simplicity, we  reduce the study to the case $\alpha\neq 2$ which breaks the translation symmetry so that we do not get involved with  equivariant bifurcation theory.
 
The Jacobian at $\mv= 0$ is given by (using some symmetries):
\[
-{\rm Id}_{3\times 3}+\lambda s_1 \mK,
\]
where $s_1 =S_0^{(1)}(0)=\frac{1}{4}$ and the matrix $\mK$ is equal to
\[
 \mK=\left[ 
 \begin{array}{ccc}
  \varepsilon_0  & \varepsilon_0\sqrt {|J_1|}\left\langle \text{1},\cos_\alpha\right\rangle&0\\
  
  \varepsilon_1\sqrt {|J_1|}\left\langle\text{1},\cos_\alpha\right\rangle &    J_1\left\langle \text{1},\cos^2_\alpha\right\rangle&0\\
  
   0 &0&  J_1\left\langle \text{1},\sin^2_\alpha\right\rangle
 \end{array}\right]
\]
$\mK$ has in general (for $\alpha\approx 2$) three real eigenvalues. Indeed the operator $J$ is self-adjoint for the inner-product defined above, but, as previously mentioned, the basis $(X_0,X_1,X_2)$ is not orthogonal for this inner product, and hence the matrix $\mK$ is not symmetric in this basis. We note $\sigma_1$ the eigenvalue of $\mK$ corresponding to the eigenvector $(0,0,1)$, $\sigma_2$ and $\sigma_3$ the two eigenvalues ot its upper lefthand $2 \times 2$ submatrix. The values, noted $\lambda_i$, $i=1,2,3$, corresponding to potential  \footnote{The upcoming nonlinear analysis will show that they are indeed bifurcation points.} bifurcations are equal to $4/\sigma_i$.
The signs of the $\lambda_i$s  give the number of bifurcated branches (recall that $\lambda>0$). Because $s_2=S^{(2)}_0(0)=0$ and $s_3=S^{(3)}_0(0)=-1/8\neq 0$, all branches are Pitchfork branches (see section.\ref{section:simpler}) whose third order term is $\chi_3^{(i)}=\lambda_i^2\frac{s_3}{6s_1}\left\langle e_i^3,e_i^*\right\rangle_2$($\approx\lambda_i^2\frac{s_3}{6s_1}\|e_i^2\|_2^2<0$ for $\alpha\approx 2$). Hence these branches are directed toward $\lambda>\lambda_i$.
This is summarized in table \ref{table:bifring}:
\begin{table}
\begin{center}
\begin{tabular}{|l|c|c|}\hline
\backslashbox{$\varepsilon_0$}{$\varepsilon_1$} &-1&1\\\hline
$-1$ & 0 & 2\\\hline
$1$ & 1 & 3\\\hline
\end{tabular}
\end{center}
\caption{Number of bifurcated Pitchfork branches from $(0,\lambda)$ depending on the values of $\varepsilon_0,\,\varepsilon_1$. The value of $\alpha$ in \eqref{eq:RM3d} is close to $2$. }
\label{table:bifring}
\end{table}
The eigenvectors $e_i$, $i=1,2,3$, of the Jacobian of \eqref{eq:RMPG1} at $\mv=0$ ar given by: 
\[
 e_1 = \sin_\alpha, \quad
 e_{2,3} = a_{2,3}+b_{2,3}\cos_\alpha
\]

We reduce the number of possibilities by assuming from now on that $J_1>0$. It turns out that in this case $\sigma_2<0$ and there are only two possibilities to consider: $\sigma_1<\sigma_3 $ for $\alpha>2$
  and $\sigma_3<\sigma_1 $ for $\alpha<2$. This gives the relative position of the different bifurcated branches, noted $P_i$, $i=1,3$.
Once we have found the bifurcation point, we can numerically compute the bifurcated branches for all positive values of $\lambda$ using a continuation method (we used the pseudo-arc length method as described for example in \cite{sala-heroux-etal:04,kuznetsov:98}). 

From our numerical experiments we conjecture that in the case $\varepsilon_0=-1$ and $\varepsilon_1=1$, $P_1$ and $P_3$ satisfy the following properties
\begin{enumerate}
	\item $P_1$ lies on the $v_3$-axis.
 	\item $P_3$ lies in the plane of equation $v_3=0$.
	\item $P_1$ and $P_3$ do not intersect.
\end{enumerate}

\begin{remark}
We can reverse the orientation of the pitchforks by choosing a nonlinearity such that $s_3>0$, the bifurcation diagram would be more complex: more saddle node points would appear because of proposition \ref{prop:saddle}. It is the fact that $s_2=0$ for the sigmoid which produces pichfork branches. Another choice of nonlinearity, for example $S(x)=\frac{1}{1+exp(-x+\epsilon)}$, would produce transcritical branches. It is indeed difficult to imagine that $s_2=0$ for an experimental fit of a ``real'' rate function. But anyway we are not yet looking at the ``real'' bifurcation diagram since we are assuming $\varepsilon=\mu=0$.
\end{remark}

\begin{figure}[htbp]
\centerline{
         \includegraphics[width=0.5\textwidth]{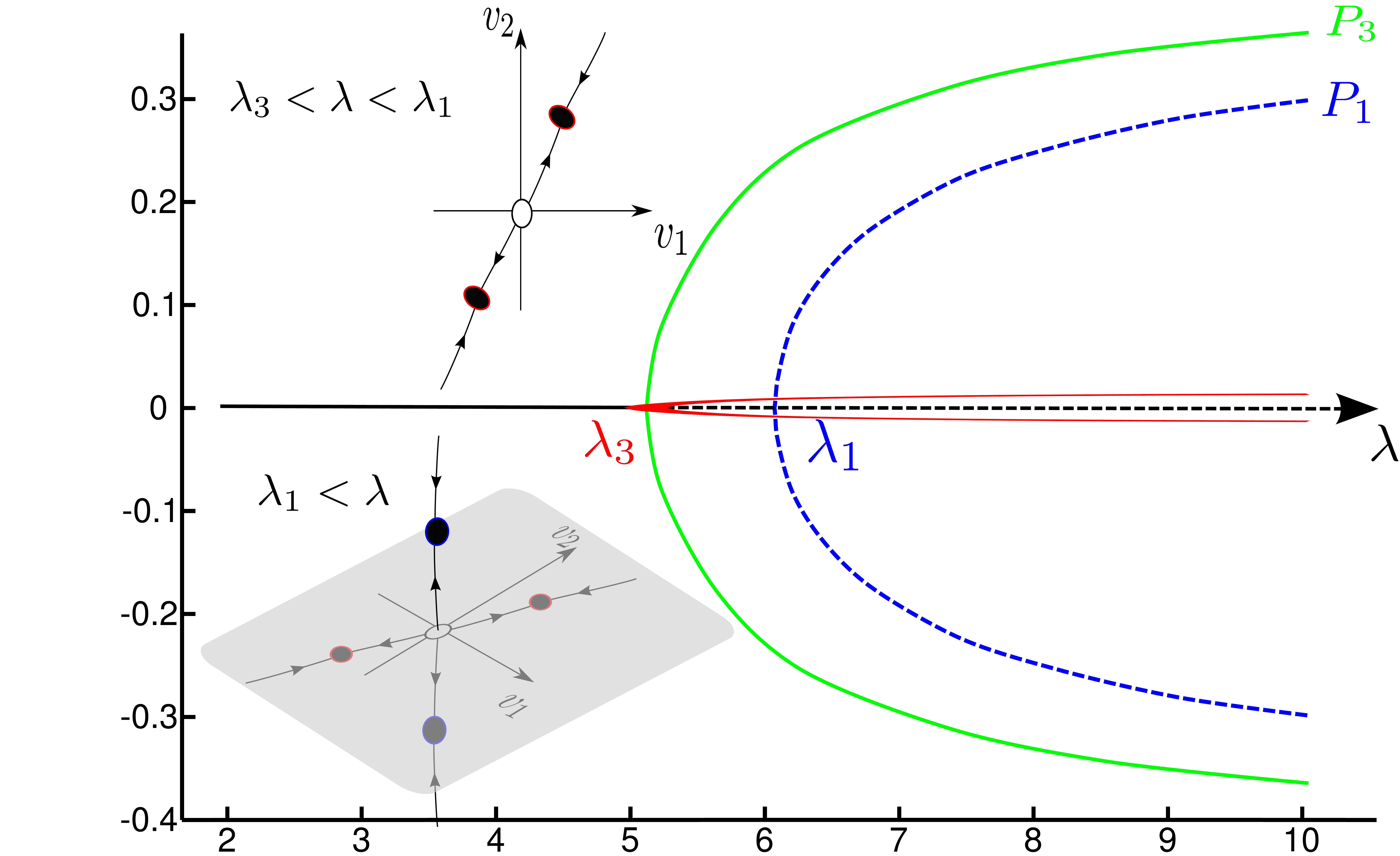}
         \includegraphics[width=0.5\textwidth]{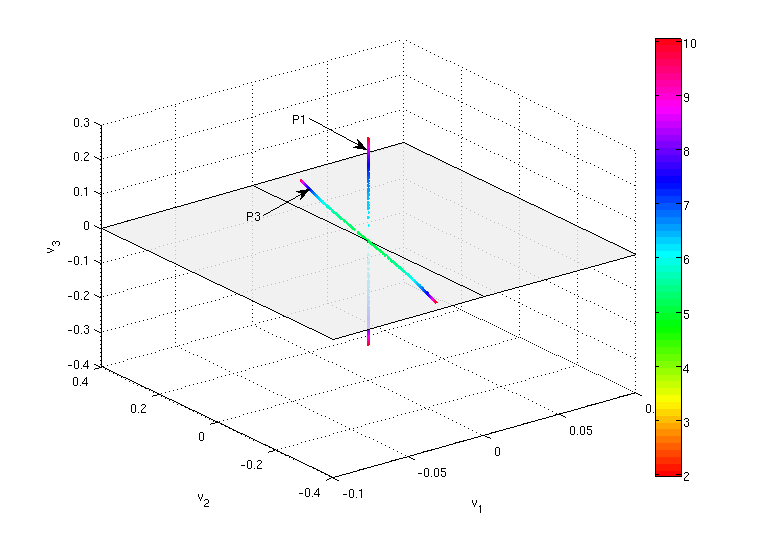}
}
       \caption{Left: Plot of the bifurcation diagram for $\alpha=2.2$. It shows the two pitchfork branches $P_1$ and $P_3$. For each branch, we have only plotted the nonzero coordinates with $v_1=red,\ v_2=green$ for $P_3$ and $v_3=blue$ for $P_1$. We have also plotted the dynamics in two and three dimensions according to the values of  the slope parameter $\lambda$. Right: Plot of the equilibrium points. The color encodes the value of the slope $\lambda$ (see text). $J_0=-1,\ J_1=1.5,\ \mu=0,\ \varepsilon=0,\ \alpha=2.2$}
\label{fig:RMA2p23D}
\end{figure}

Figure \ref{fig:RMA2p23D} shows a typical example corresponding to the values of the parameters that are found in the largest number of published articles.  We  come back to this choice in section \ref{section:RMDiscussion}. The left part of the figure shows the three components of the persistent states as functions of $\lambda$. For $P_1$ there is only one nonzero component, $v_3$, in blue. For $P_3$ there are two nonzero components, $v_1$ shown in red and $v_2$ shown in green. The right part of the figure shows another representation of $P_1$ and $P_3$ as curves parametrized by $\lambda$ in the $(v_1,v_2,v_3)$ space. $P_3$ is clearly in the $(v_1,v_2)$ plane while $P_1$ is along the $v_3$-axis. The color at each point of the curves represents the value of $\lambda$ according to the color scale shown on the right.

In detail we have
\begin{description}
 	\item[$\lambda_{3} < \lambda < \lambda_1$] When $\lambda$ goes through $\lambda_3$, the 0-solution loses its stability and becomes a saddle. There are three persistent states, $0$ (unstable node) and two located on the pichfork branch $P_3$, both stable. The corresponding dynamics is shown in the left part of figure \ref{fig:RMA2p23D}.
 	\item[$\lambda_{1} < \lambda$] When $\lambda$ goes through $\lambda_1$, the 0-solution loses its stability along the $v_3$-axis. There are two new persistent states located  on the pitchfork branch $P_1$, both are unstable nodes (the unstable manifold is one-dimensional). The corresponding dynamics is also shown in the left part of figure \ref{fig:RMA2p23D}.
\end{description}
There are at most $5$ stationary solutions.

%

\begin{remark}
 In the case $J_0=1$ and $\varepsilon_1=1$, there is another pichfork branch to handle, see table \ref{table:bifring}. For $\lambda$ big enough, $\mv=0$ becomes an unstable node and there are $7$ stationary solutions instead of 5 in the case $J_0=-1$
\end{remark}

\subsubsection{The case $\mu=1,\ \varepsilon\neq 0$}
We are now halfway from our scheme completion. To have an idea of the persistent states at low contrast (\textit{ie} $\varepsilon\approx0$), we need to know the persistent states for : $\lambda,\mu=1,\varepsilon=0$, that is we need to know the solutions of 
\[
 V^f_\lambda=J \cdot S(\lambda V^f_\lambda)+\theta
\]

Following our program, we numerically compute the persistent states when the slope $\lambda$ and $\mu$ both vary. As described in section \ref{subsection:blambda} we expect many of the previous bifurcations to disappear thereby breaking some of the connectivity of the sets ${\mathcal B}_\lambda$ which can actually be (partially) recovered by considering the sets ${\mathcal B}_{\lambda,\mu}$
This was done using the library TRILINOS (see \cite{sala-heroux-etal:04} and the \href{http://trilinos.sandia.gov}{website}) using multiparameters continuation. 

\begin{figure}[htbp]
\centerline{
	\includegraphics[width=0.5\textwidth]{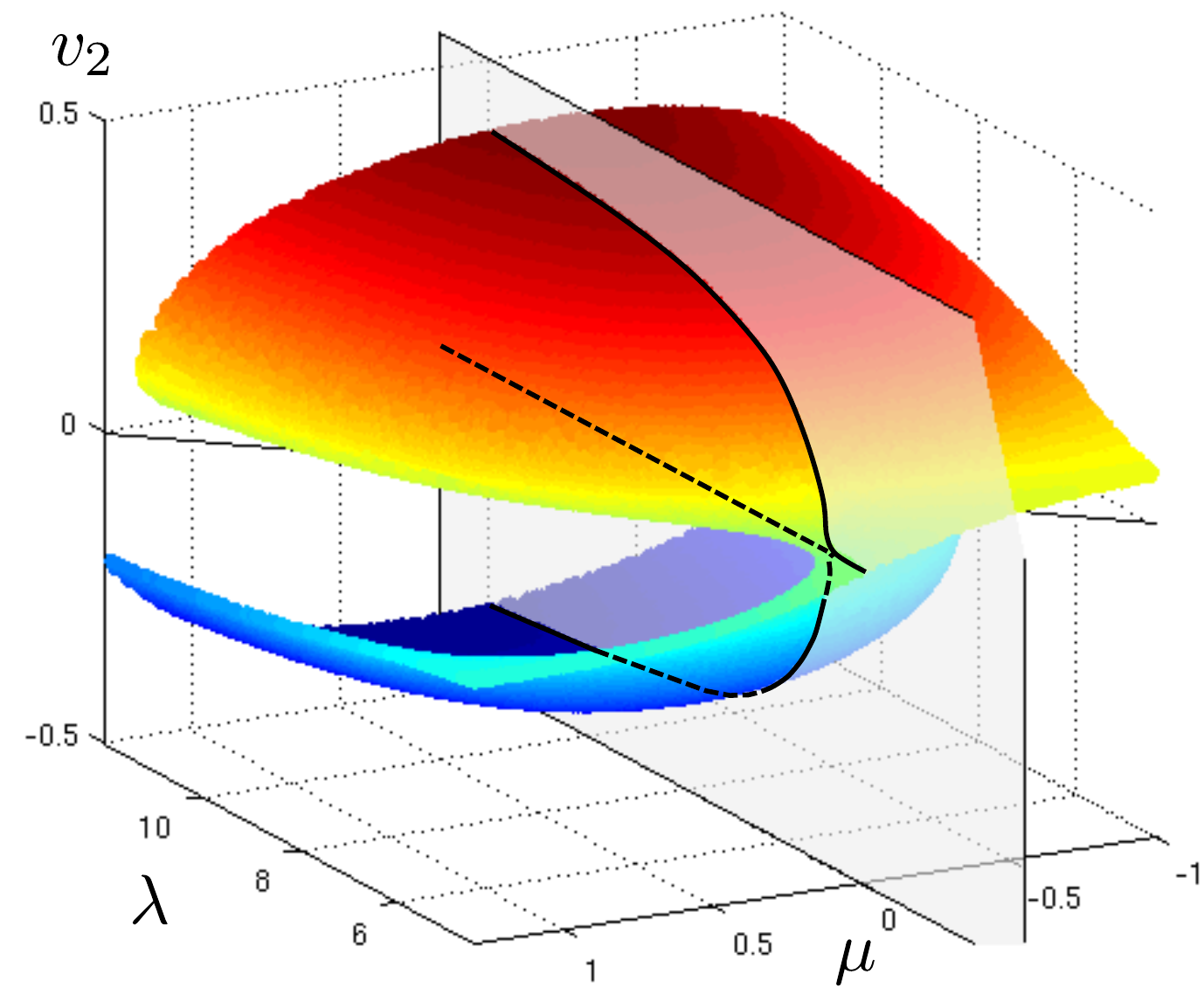}
        \includegraphics[width=0.5\textwidth]{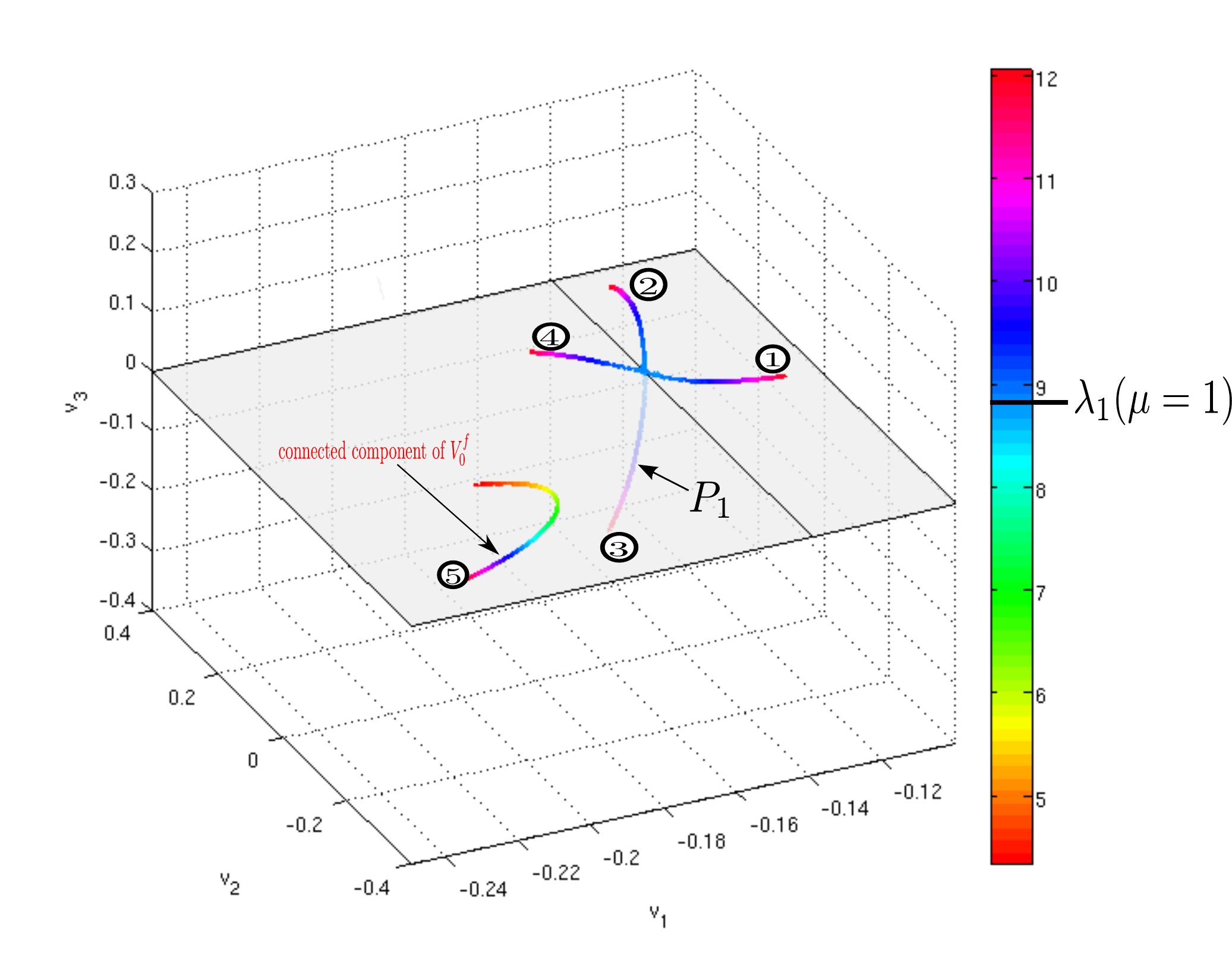} 
}
       \caption{Left: The $v_2$ component of the 2-parameters continuation $(\lambda,\mu)$ for $\theta=0.1,\ \alpha=2.2,\ J_0=-1,\ J_1=1.5$. Right: Plot of the equilibrium points. The color encodes the value of the slope $\lambda$. $J_0=-1,\ J_1=1.5,\ \mu=1,\ \varepsilon=0,\ \alpha=2.2,\ \theta=0.1$. }
       \label{fig:RMmu1}
\end{figure}

We show an example of this continuation in figure \ref{fig:RMmu1}.Left where we display the $v_2$ component of the persistent states as a function (sometimes multivalued) of $\lambda$ and $\mu$. A cross-section of this set by the plane of equation $\mu =0$ (shown as semi-transparent in the figure) yields a curve identical to the one shown in green figure \ref{fig:RMA2p23D}.Left. The figure nicely shows how the first pichfork bifurcation branch $P_3$ opens up when $\mu$ becomes non zero: this gives the connected component of $\vf$ which is linearly stable.
  
\begin{equation}
 \begin{cases}
 \label{eq:RMMU}
   0=v_1- \varepsilon_0 \left\langle S_0(\lambda V),\text{1}\right\rangle+\mu(\theta+\frac{\varepsilon_0}{2})\\
   0=v_2 -   \varepsilon_1\sqrt{|J_1|}\left\langle S_0(\lambda V),\cos_\alpha\right\rangle+\mu\varepsilon_1\sqrt{|J_1|}\frac{sin_\alpha(\pi/2)}{\alpha\pi}\\
   0=v_3 - \varepsilon_1\sqrt{|J_1|}\left\langle S_0(\lambda V),\sin_\alpha\right\rangle
  \end{cases} 
\end{equation}

However, as can be seen from \ref{fig:RMmu1}.Right, non-zero values of $\mu$ do not break the pitchfork $P_1$. It is easy to qualitatively understand why, even though a full mathematical proof is hard to come up with: $\mu$ does not affect the third equation which produces the pichfork $P_1$. We can prove it locally for $\mu$ near $0$ using the implicit function theorem. We are looking for a point $(v_1(\mu),v_2(\mu),0)$ at which a pichfork occurs for $\lambda=\lambda_1(\mu)$. Let us consider 
\begin{equation}
 H(v_1,v_2,\lambda;\mu)=\left[ 
 \begin{array}{l}
  v_1- \varepsilon_0 \left\langle S_0(\lambda V),\text{1}\right\rangle+\mu(\theta+\frac{\varepsilon_0}{2})\\
  v_2 -   \varepsilon_1\sqrt{J_1}\left\langle S_0(\lambda V),\cos_\alpha\right\rangle+\mu\varepsilon_1\sqrt{|J_1|}\frac{sin_\alpha(\pi/2)}{\alpha\pi}\\
  1-\lambda J_1\left\langle DS_0(\lambda v_1X_0+\lambda v_2 X_1),sin_\alpha^2\right\rangle
 \end{array}
\right] 
\end{equation}
 where the last component of $H(v_1,v_2,\lambda;\mu)$ is $\frac{\partial}{\partial v_3}(\ref{eq:RMMU}.3)$. It is easy to see that $H(0,0,\lambda_1,0)=[0\ 0\ 0]$. The Jacobian of $H$ w.r.t. $(v_1,v_2,\lambda)$ at $(0,0,\lambda_1,0)$ is (because $S^{(2)}(0)=0$) : 
\[
 \left[ 
 \begin{array}{ccc}
  -1+\lambda_1s_1\varepsilon_0  & \lambda_1s_1\varepsilon_1\sqrt {|J_1|}\left\langle \text{1},\cos_\alpha\right\rangle&0\\
  
  \lambda_1s_1\varepsilon_0\sqrt {|J_1|}\left\langle\text{1},\cos_\alpha\right\rangle &    -1+\lambda_1s_1J_1\left\langle \text{1},\cos^2_\alpha\right\rangle&0\\
  
   0 &0&  -s_1J_1\left\langle \text{1},\sin^2_\alpha\right\rangle
 \end{array}\right]\in GL_3(\mathbb R)
\]
Hence there exists a unique solution defined locally for $\mu\geq 0$ satisfying $H(v_1(\mu),v_2(\mu),\lambda_1(\mu);\mu)=0$ : we have found  a bifurcated point. Moreover, as $\chi^{(3)}_1(\mu=0)\neq0$, it will remains so for small $\mu$: the Pichfork $P_1$ is not affected by $\mu$. For large values of $\mu$ we have to rely on numerical simulations.

\begin{figure}[htbp]
\centerline{
	\includegraphics[width=0.5\textwidth] {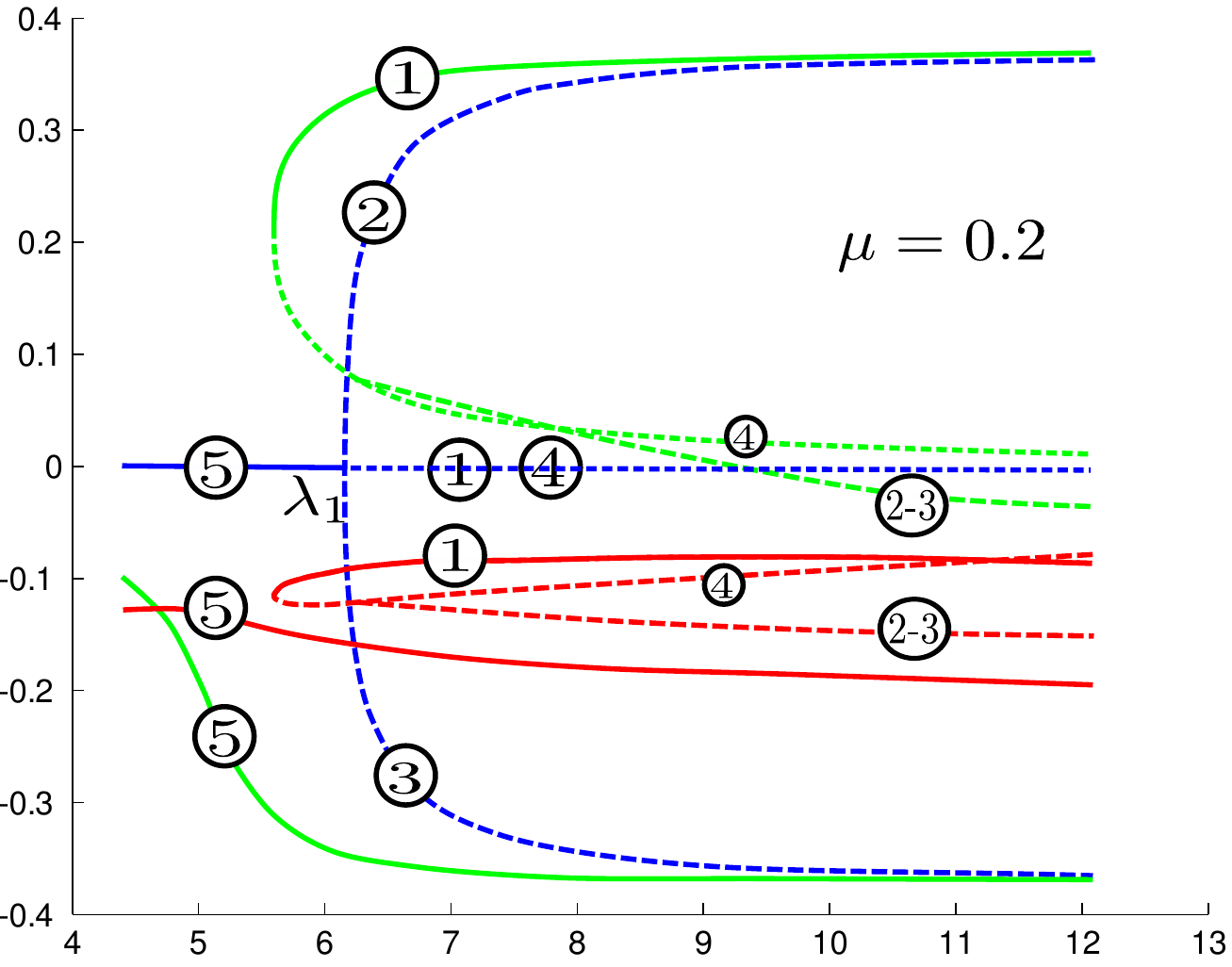}
	\includegraphics[width=0.5\textwidth]{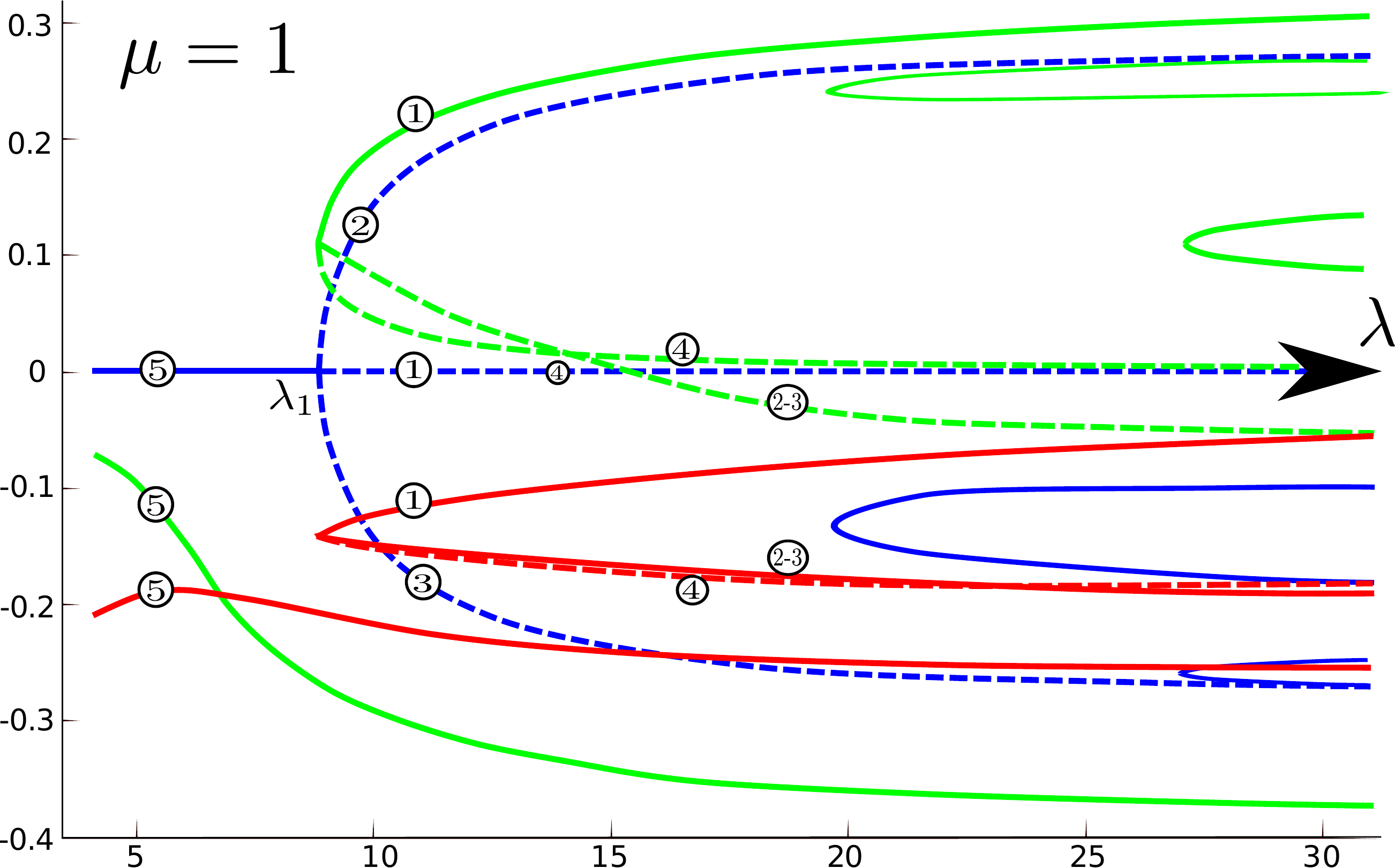}
}
       \caption{ Plot of the three components $(v_1,v_2,v_3)=(red,green,blue)$ as functions of the slope $\lambda$ for the following values of the parameters: $J_0=-1,\ J_1=1.5,\ \mu=1,\ \varepsilon=0,\ \alpha=2.2,\ \theta=0.1$. Recall that the slope of the model has to be around $\lambda=29$.}
       \label{fig:RMComponentsMu1}
\end{figure}

Now, because of lemma \ref{lemma:sym}, the solutions $(v_1,v_2)$ corresponding to $v_3\neq0$ (\textit{ie} lying on the pichfork branch) are the same for $v_3$ and $-v_3$ which gives the branch 2-3 in figure \ref{fig:RMComponentsMu1}. Hence, when $\mu\neq 0$, we still have the pichfork $P_1$ and the branch (located in $v_3=0$) arising from the opening of $P_3$. This gives the  diagram shown in figure \ref{fig:RMComponentsMu1} for the three components of $\mv$.

This diagram is a bit misleading because if we count the $red$ components, there are four of them for $\lambda>\lambda_1$ which gives an even number of solutions (in contradiction with proposition \ref{prop:bumps}). In fact by doing so we miss the symmetry $v_3\leftrightarrow-v_3$ and the corresponding solutions. It is easier to look at figure \ref{fig:RMmu1}.Right to count the stationary solutions.
\

From section \ref{section:simpler}, it follows that the connected branch of $V_0^f$ is stable as well as the branch $P_1$. The only unstable branch comes from the opening of $P_3$ and is shown in figure \ref{fig:RMmu1}.Right.

\begin{remark}
 Figure \ref{fig:RMmu1} gives a good example where our scheme allows to detect another branch of solutions which is not connected to the connected component of the trivial solution.
\end{remark}

\begin{figure}[htbp]
\centering
         \includegraphics[width=\textwidth]{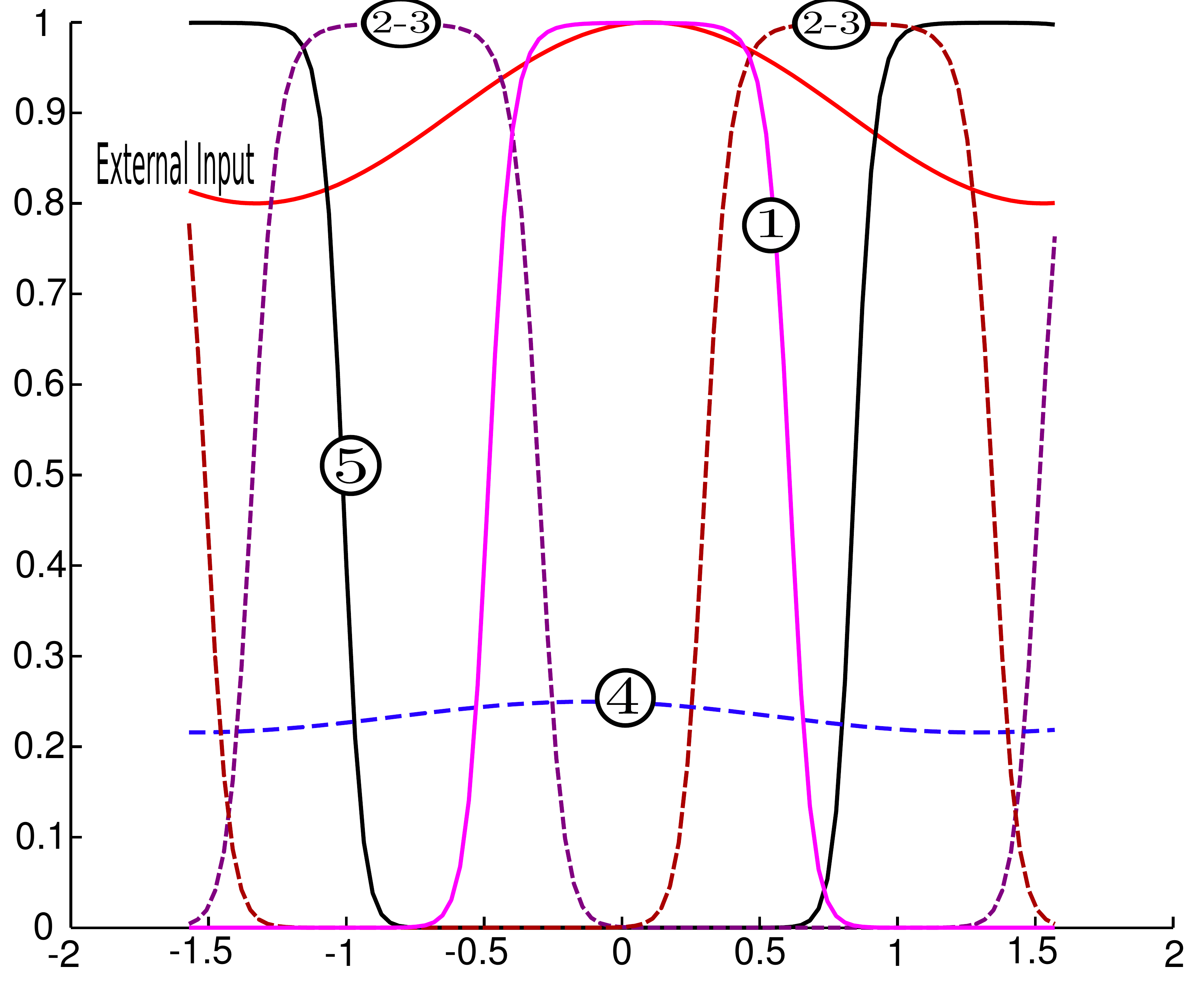}
       \caption{Grey: Stationary solutions $\mA^f$ for $\varepsilon I(x)=\varepsilon(1-0.1+0.1cos_\alpha(x-0.1)),\ \alpha=2.2.\ J_0=-1,\ J_1=1.5,\ \varepsilon=0.1, \lambda=29.2$. Red: plot of $\mI$. Blue: solution in the connected component of $\mV^f_0$.}
\label{fig:RMLambda9}
\end{figure}

The figure \ref{fig:RMComponentsMu1} also tells us which branches will appear when the contrast satisfies $\varepsilon\neq 0$: this will be a perturbation of figure \ref{fig:RMmu1}. Moreover, our scheme detects two saddle-node at $\lambda=20$ and $\lambda=26$ which are undetectable by $\lambda-$continuation only. Each saddle-node gives 2 additional persistent states : one stable, the other unstable.
Except in the case $x_0=0$ (see equation \eqref{eq:RM3d}), the pichfork $P_1$ will open up, giving two other connected components in addition to the connected component of $\mV_0^f$. As an example, we have plotted in figure \ref{fig:RMLambda9} all the persistent states (in the activity representation, \textit{ie} $\mA^f=S(\mV^f)$, in the case $I=1-0.1+0.1+cos_\alpha$ and $\varepsilon=0.1,\ \lambda=9$. There are five of them as predicted from a perturbation of figure \ref{fig:RMComponentsMu1}.

\subsection{Discussion}\label{section:RMDiscussion}
There are two reasons why we presented this example. First, it is a nice simple model to which the formalism of this article easily applies and allows us to push the analysis far enough to grasp an almost complete understanding of its persistent states and a somewhat detailed understanding of its dynamics. Second, it conveys information for models of V1 that is likely to be biologically relevant. For example, as the slope $\lambda$ of the sigmoid is increased, many (up to 4 depending on the signs of $(J_0,J_1)$) new stationary states appear whose stability evolves with $\lambda$. One of these solutions is ``dramatic'' for the purpose of orientation detection: even if the LGN input orientation peaks around the angle $x_0=0$, the Ring Model can produce a stable cortical state (or a percept) corresponding to an angle of $\pi/2$! 

However these solutions may be destabilized by adding lateral spatial connections in a spatially organized network of ring models; it remains an area of future investigation. As far as we know, only Bressloff and co-workers looked at this problem  (see \cite{bressloff-cowan:02c,bressloff-cowan-etal:01}): they studied the local dynamics around the bifurcation points but did not look at non-local dynamics (basically they considered that there were 2 stationary solutions around a bifurcation point and forgot about the third one predicted by proposition \ref{prop:bumps}). 

There is no  biologically motivated restriction on the values of the slope $\lambda$ which can be as large as desired, without mentioning the fact that neural mass models are very often  written with a Heaviside function for the nonlinearity which, as mentioned previously, is the limit case: $H(x)=\lim\limits_{\lambda\to\infty}S(\lambda x)$.

We also made the assumption $\alpha\neq 2$ in the previous analysis. It remains to know how much of what preceed holds in the case $\alpha=2$. More generally, what would remain if one were to choose a difference of Gaussians as a connectivity function over a cortex $\Omega=(-\pi/2,\pi/2)$. We expect a lot of similarities if the width of the difference of Gaussians is of the same order of $cos_\alpha$.

\section{Two populations of spatially organized neurons}\label{section:numerics2d}
We now apply the previous analysis to a system we started to analyze in \cite{faugeras-veltz-etal:09}. In the somewhat reduced form we consider here, it consists of two populations ($p=2$), one excitatory, one inhibitory, distributed over a flat ($d=2$) cortex $\Omega=[-1,\ 1]^2$.  The connectivity matrix kernel writes:
\begin{equation}
\label{eq:2pop2d}
\mJ(\mr,\mr')=\left[ 
\begin{array}{cc}
 a G_{11}(\mr-\mr')&-bG_{12}(\mr-\mr')\\
 bG_{21}(\mr-\mr')&-cG_{22}(\mr-\mr')
\end{array}
\right]  
\end{equation}
where the $G_{ij}(\mr)=e^{-\frac{\|\mr\|^2}{2\sigma_{ij}}}$ are two-dimensional Gaussian functions defined on $\mathbb R^2$ with $G_{12}=G_{21}$. $a,\,b,\,c>0$ characterize the strength of the connections. We also assume $\mI_{ext}=0$. The parameter $\mu$ controlling the translation of the sigmoid (see \eqref{eq:net3}) is therefore the only parameter, outside $\lambda$, that we will vary from 0 to 1. We also chose (notice that $s_2\neq 0$) :
\[
 S(x)=\frac{1}{1+e^{-x+\theta}},\ \ \theta=1.3
\]

In \cite{faugeras-veltz-etal:09}, we were able to compute the stationary solutions $\mV^f_{\lambda,\mu=1}$ when the slope $\lambda$ was small (\textit{i.e.} $\lambda<\lambda^*$) using the Nystrom method. We now know that if we perform a continuation of these solutions with respect to $\lambda$, we are bound to miss quite a few of them. Therefore we perform a two-parameter continuation with respect to the pair $(\lambda,\mu)$ in order to recover more, if not all, stationary solutions.  

Biologically speaking, no systematic investigation has been performed to test the validity of the translation invariance of the connectivity function. Hence a roughly translation invariant (called heterogeneous in \cite{ermentrout-cowan:80}) is not less biologically relevant. This is where the PG-kernels are useful: they provide an easy way to approximate the convolution operation as well as an effective representation of the connectivity (see section \ref{section:PG}).

There are four reasons why we think this example is interesting: 
\begin{itemize}
 \item We want to show how to deal with heterogeneous kernels.
 \item We want to give a non trivial example to the existence of a branch of solutions not connected to the trivial solution $\mV^f_0$.
 \item We want to show an example of application of proposition \ref{prop:saddle}.
 \item More importlantly, we want to show how the results of section \ref{section:rabinowitz} may change in the two-dimensional case.
\end{itemize}

As in the previous case, we are not interested in having the complete bifurcation diagram of the system but rather in giving numerical examples of the previously enumerated points.

\subsection{PG-kernel approximation of $\mJ$}
Let us write
\[
e^{-\|\mr-\mr'\|^2/2}=e^{-\|\mr\|^2/2}e^{-\|\mr'\|^2/2}e^{\langle\mr,\mr'\rangle}\approx e^{-\|\mr\|^2/2}e^{-\|\mr'\|^2/2}\Big(1+\langle\mr,\mr'\rangle+\frac{1}{2}\langle\mr,\mr'\rangle^2+... \big)
\]
We notice two important facts: 
\begin{itemize}
 \item $1+\langle\mr,\mr'\rangle+\frac{1}{2}\langle\mr,\mr'\rangle^2+...$ is a polynomial in the components of $\mr$ and $\mr'$, hence a PG-kernel.
\item $e^{-\|\mr\|^2/2}$ has a bell-shape which we approximate with the following function $\Big(1-\|\mr\|^2\Big)^a$ with $a>0$. This choice is also motivated by the fact that $e^{-\|\mr\|^2/2}$ tends to zero at the edges of the infinite cortex that is usually considered in the literature. We keep this property with our finite cortex since $\Big(1-\|\mr\|^2\Big)^a = 0$ on  the boundary $\partial \left([-1,1]^2\right)$ of $[-1,1]^2$. Any other positive bell shaped function would be appropriate.
\end{itemize}
Putting these two facts together, we end-up  with the following approximation of a Gaussian convolution kernel with a PG-kernel.
\begin{equation}
\label{eq:PG}
J_{ij}(\mr ,\mr')=C_{ij}+\Big(1-\|\mr\|^2\Big)^{a_{ij}}\Big(1-\|\mr'\|^2\Big)^{a_{ij}'} P_{ij}(\mr,\mr'),\,i,\,j=1,\,2
\end{equation}
where $\mathbf{P}$ is a $2 \times 2$ matrix with polynomial entries in $\mr$ and $\mr'$ and $\mC$ is a constant $2 \times 2$  matrix.

\

Recall from section \ref{subsection:functional} that for $d=2$, we require $m=1$, i.e. $\mF=\mW^{1,2}(\Omega,\R^2)$. Thus $\| \mJ \|_\mF<\infty$ imposes $2(a_{ij}-1)>-1$ hence $a_{ij}>1/2$, similarly for $a_{ij}'$.

\subsection{Numerical experiments}

In these numerical experiment we choose $a_{11}=a_{11}'=3$, $a_{12}=a_{12}'=a_{21}=a_{21}'=2$, $a_{22}=a_{22}'=4$, $a=10,b=15,c=12.75$ and $\mL=0.5Id$. The polynomials $P_{ij}(\mr,\mr')$ are equal to the 4th-order Taylor expansion of $e^{\frac{\langle\mr,\mr'\rangle}{\sigma_{ij}}}$, with $\sigma_{ij}^{-1}=a_{ij}$.

This results in the following expressions for the two components $V^f_1$ and $V^f_2$ of the persistent state $\mV^f$:
\[
 \begin{array}{c}
  V^f_1=Q_{11}(\mr)\Big(1-\|\mr\|^2\Big)^3+Q_{12}(\mr)\Big(1-\|\mr\|^2\Big)^2\\
V^f_2=Q_{21}(\mr)\Big(1-\|\mr\|^2\Big)^2+Q_{22}(\mr)\Big(1-\|\mr\|^2\Big)^4
 \end{array}
\]
The polynomials $Q_{ij}$, $i,\,j=1,2$ are of degree 4. The total number of parameters needed for representing a persistent state is therefore  equal to $5\times5\times4=100$ variables. This number comes from the fact that we have four polynomials $Q_{ij}$, $i,j=1,2$, of total degree 10 in two variables that are the products of two polynomials in one variable of degree 4. The reason for this is the special form of the polynomials $P_{ij}(\mr,\mr')$ which are functions of the inner product $\langle \mr,\mr' \rangle$. This is a somewhat crude approximation of the convolution kernel (\ref{eq:2pop2d}).

\begin{figure}[htbp]
\centering
         \includegraphics[width=\textwidth]{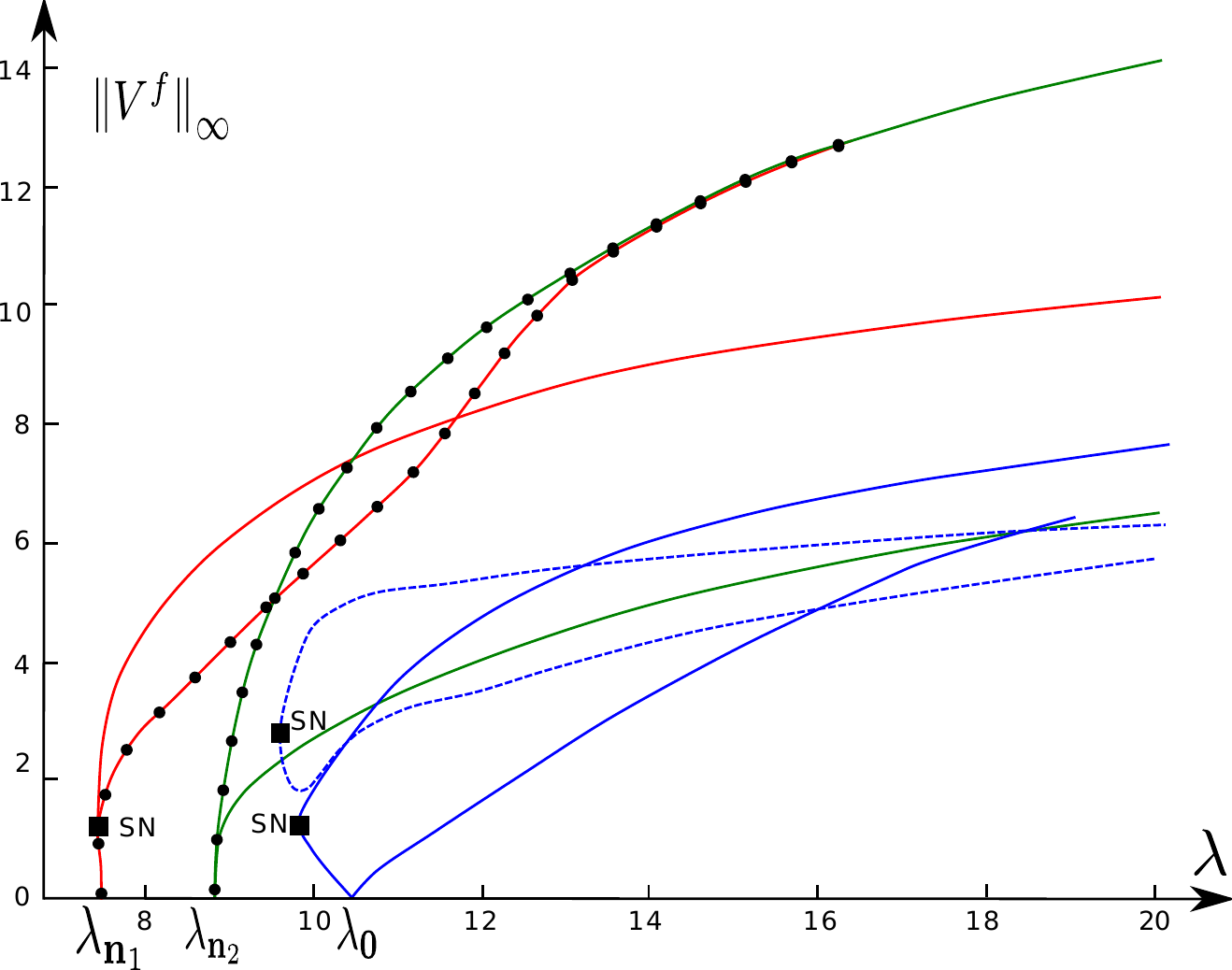}
       \caption{Plot of the infinity norm $\|\cdot\|_\infty$ of the 100-dimensional vector $\mV^f$ as a function of the slope parameter $\lambda$ for $\mu=0$. Notice that $\mV^f=0$ is solution for all values of $\lambda$. The dotted branch (red and green) connects $\lambda_{\mathbf{n}_1}$ to $\lambda_{\mathbf{n}_2}$. The dashed-line (blue) is a set of solutions that is not connected to the component corresponding to the trivial solution $\mV^f=0$.  SN stands for saddle-node.}
\label{fig:2popinfty}
\end{figure}

Each of these 100 variables is a function of the slope parameter $\lambda$ and contrast parameter $\mu$. We represent in figure \ref{fig:2popinfty} the infinity norm of the 100-dimensional vector $\mV^f$ as a function of the slope parameter $\lambda$ for $\mu=0$.

\

All our analytical predictions are performed, for convenience and simplicity, on the convolution kernel $\mJ$. They are likely to carry over to the heterogeneous approximate case.

 The general analysis having been conducted (see section \ref{section:simpler}), we wish to do some specific computations to explain the numerical results : the first thing to do is to find the eigen-elements of $\mJ$.

For each $\mn=(n_1,n_2) \in \mathbb{N}$ consider the function $\cos_\mn : \R^2 \to \R$ defined by $\cos_\mn(\mr)=\cos(n_1 r_1+n_2 r_2)$, and the vector subspace $X_{\cos_\mn}$ of $\mF$ defined by $X_{\cos_\mn}=\{\mV \, | \, \mV= \cos_\mn \mv,\, \mv \in \R^2\}$. $X_{\sin_\mn}$ is defined in a similar fashion.

It is easy to check that the vector subspaces $X_{\cos_\mn}$ of $\mF$ are invariant by the linear operator $\mJ$, and, by parity, that $\mJ \cdot X_{\sin_\mn} =0$ for all $\mn \in \mathbb N^2$. We have
\[
 \mJ \cdot (\cos_{\mn} \mv)
 =\cos_{\mn}\ \hat\mJ(\mn) \mv=
\cos_{\mn} \left[ 
\begin{array}{cc}
 a \hat G_1(\mn)&-b\hat G_2(\mn)\\
 b\hat G_2(\mn)&-c\hat G_3(\mn)
\end{array}
\right]  \mv \quad \forall \mv \in \R^2
\]
where $\hat \mJ$, $\hat G_i$ are the Fourier transforms of $\mJ$ and $G_i$. 

Keeping the notations of section\ref{section:simpler}, we find that 
\[
 \sigma_\mn\in\bigcup_{\mm \in\mathbb N^2}\quad\sigma\Big(
\hat\mJ(\mm)\Big)
\]
and the corresponding eigenvectors\footnote{We do not prove that we obtain all of them though it appears that these particular eigen-elements are sufficient to explain the numerically observed facts.} are in $X_{\cos_\mn}$ for some $\mn \in \mathbb N^2$, hence of the form $\cos_\mn \mv$ for some $\mv \in \R^2$.
We have chosen the Gaussians $ G_i$ such that the first eigenvalues are simple, hence each point $\lambda_\mn$ is a bifurcation point: the corresponding \textit{chi}-factors are 
\[
 \chi_2^{(\mn)}\propto\left\langle \cos_{\mn}^2,\cos_{\mn}\right\rangle =
  0\text{ if } \mathbf{n}\neq\vec 0\\
 \]
 \[
  \chi_2^{(\mathbf{0})}\neq 0
 \]
Hence $\lambda_{\mathbf{0}}$ is a transcritical bifurcation point whereas the other points $\lambda_{\mathbf{n}}$ are pitchfork bifurcation points. Depending on the vector $\mv$, the orientation of the pichfork may change.

In figure \ref{fig:2popinfty} (case $\mu=0$), we see two pichforks at $\lambda_{\mn_1},\ \lambda_{\mn_2}$ and the transcritical branch at $\lambda_\mathbf{0}$. The first pichfork at $\lambda_{\mn_1}$ (red branch) is oriented toward the decreasing $\lambda$s, thus according to proposition \ref{prop:saddle}, a saddle-node (noted SN in the figure) must appear. The same is true for the transcritical branch (continuous blue branch). Locally (near $\lambda_\mathbf{0}$), each coordinate $V^f_i, i=1,2$ is 'flat', \textit{ie} it looks like $cos_\mathbf{0}$.

\begin{figure}[htbp]
\centering
         \includegraphics[width=\textwidth,angle=-90]{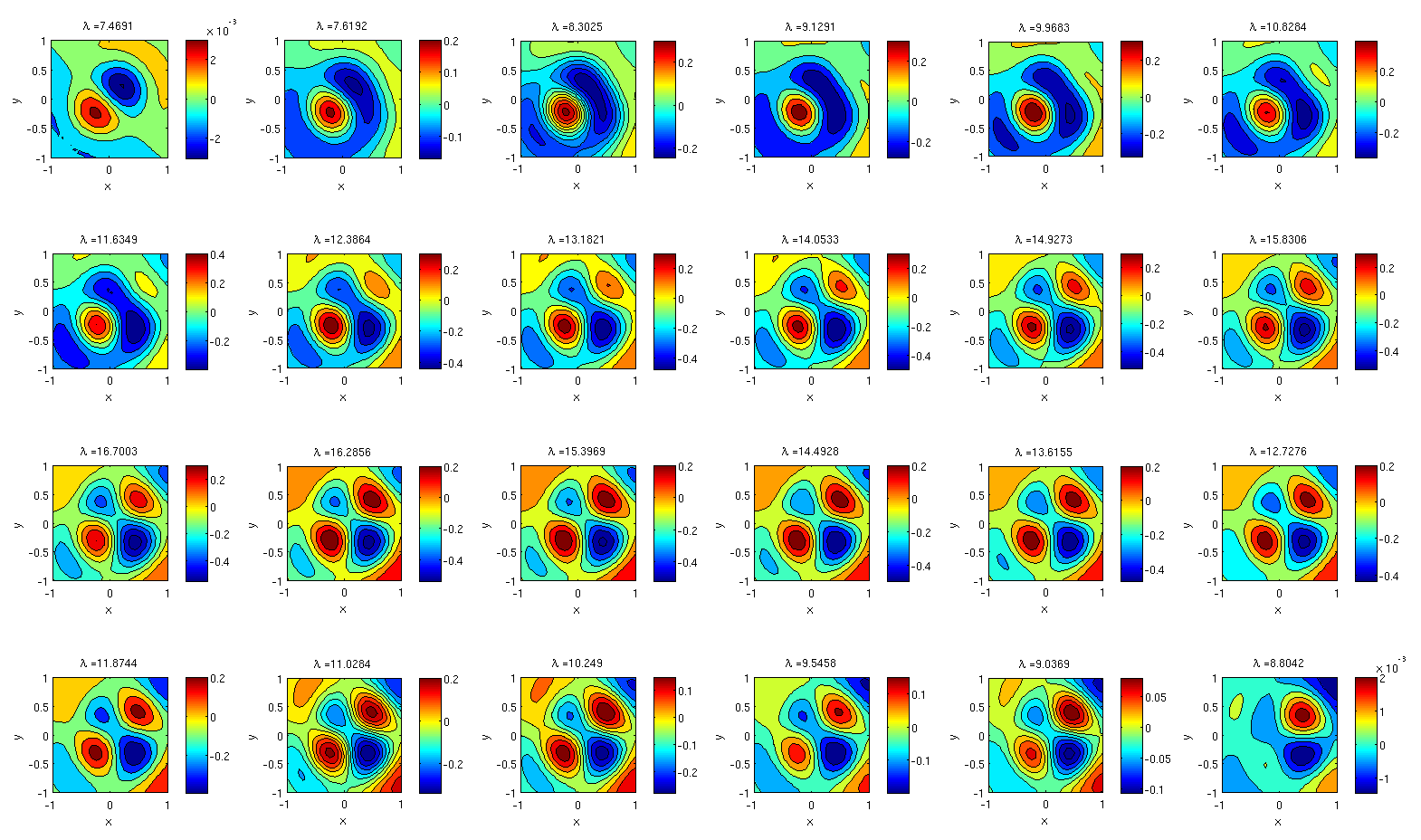}
       \caption{Plot of the first component $V^f_1$ of the solution $\mV^f$ along the red-green dotted branch connecting $\lambda_{\mn_1}$ to $\lambda_{\mn_2}$ in figure \ref{fig:2popinfty}.}
\label{fig:branche}
\end{figure}

The bifurcation points $\lambda_{\mn_1}$ and $\lambda_{\mn_2}$ are connected by the red-dotted branch. It cannot be seen from this graph that this is true (because two states with the same norm may be different), but it can be checked by looking at the 100 components. We plot in figure \ref{fig:branche} the stationary membrane potential $V^f_1$ of the first population for different values of $\lambda$ along the branch connecting $\lambda_{\mn_1}$ to $\lambda_{\mn_2}$. 

The last interesting fact is that our multi-parameter scheme allows to compute branches that are not connected to the trivial solution (here $\mV^f=0$) as one can see in figure \ref{fig:2popinfty}: the dashed blue branch cannot be detected by $\lambda$-continuation of the trivial solution $\mV^f=0$ because it does not intersect any branch coming from $\mV^f=0$. This, again, cannot be read directly from figure \ref{fig:branche}: one has to look at all 100 components.

\begin{remark}
 We have chosen the parameters $a,b,c$ such that the first three eigenvalues have zero imaginary parts. Numerically, most of the other eigenvalues $\sigma_\mn$ have non-zero imaginary part leading to Hopf bifurcations. Other stationary bifurcation may appear for large slope values (\textit{ie} $\lambda>40$).\\
\end{remark}

\section{Conclusion}
In this paper we have pursued the analysis, started in \cite{faugeras-grimbert-etal:08,faugeras-veltz-etal:09}, of a special type of integro-differential equation that appears in neural field and neural mass models where we are interested in approximating mesoscopic and macroscopic ensembles of neurons by  continuous descriptions.  We call these equations neural equations. Like in these two previous papers our approach is based on functional analysis to obtain clean results about the well-posedness of these equations from which we can study the questions of existence and uniqueness of the solutions. In this article we have enforced more spatial regularity  by choosing the Sobolev, and Hilbert, space $\mW^{m,2}$ instead of the space ${\rm \mL}^2$ that we used previously. The reason for this change is not only technical, it allows us to study the question of the bifurcation of the solutions, while ensuring that the membrane potential (or the activity) remains finite on the cortex.

In effect these equations depend upon a number of biological or experimental parameters such as the slope $\lambda$ of the nonlinearity, the connectivity matrix $\mJ$, or the input $\mI$. These parameters vary in general in neural populations because of such processes as plasticity and learning. It is therefore important to understand how the solutions of these equations vary when these parameters change. To this end, we used two theories: the degree theory and the bifurcation theory. The degree theory does not require the above spatial regularity because it describes the general behaviour of the cortical states as the parameters vary. On the other hand, the bifurcation theory describes the precise local behaviours of the cortical states as the parameters vary  and requires  spatial regularity to be able to complete the numerical computations.

We believe that a good model should exhibit bounded membrane potentials which requires to take a bounded nonlinearity (as opposed to some papers, see for example \cite{ben-yishai-bar-or-etal:95}). The degree theory yields the powerful estimate that the number of persistent states has to be odd, hence it predicts an additional persistent state in the neighborhood of pichfork/transcritical bifurcation points. This extra point is invisible to bifurcation theory, which is a local theory (This was conjectured numerically in \cite{ermentrout-cowan:80}). It may change drastically the dynamics and shows that local analysis is not sufficient for the study of the neural field equations.

We have focused on the description of the stationary solutions of the neural equations when varying the slope parameter, and tried to compute numerically the additional persistent state given by the degree theory for arbitrary external current and connectivity. As the set of stationary solutions may not be connected, we used a multiparameter continuation scheme in order to compute non-connected branches of solutions and were able to show examples of these (see sections \ref{section:ring} and \ref{section:numerics2d}). Whether these different branches of solutions intersect or are unbounded is still unknown in the general case. However, the scalar case for a one-dimensional cortex ($p=1,\ d=1,\mI_{ext}=0$) is almost completely solved (see section.\ref{section:rabinowitz}). This point has never been mentioned in the literature to our knowledge. The question of whether we have computed all the solutions by using our multiparameter scheme is unfortunately still open.

To sum up, using a bounded nonlinearity for the firing rate introduces new stationary solutions that were not predicted before. This suggests that the analysis of neural field models (of the visual system or of the memory for example) should be re-examined. This is the focus of our current efforts.

\section*{Acknowledgments}
We thank M.E. Hendersen for his valuable help in using the library Multifario (part of Trilinos). Much of this work would have proved difficult - if not impossible - without his help. We also thank D. Wasserman for his patience and constant help in programming.

\appendix
\section{Well-posedness of operators}\label{appendix:Banach}
We prove proposition \ref{prop:FtoF}.
\begin{proof}
The integral in the righthand side of \eqref{eq:W} exists because for almost all $\mr \in \Omega$ the $p^2$ elements $W_{ij}(\mr,\cdot,t)$, $i,j=1,\cdots,p$ of $\mJ$ are in ${\rm L}^2(\Omega)$ for all $t>0$ and the $p$ coordinates of $\mV(\cdot,t)$ are in ${\rm L}^2(\Omega)$ for all $t>0$.

For each $t>0$ this righthand side defines an element of $\mW^{m,2}(\Omega)$. Because of Fubini's theorem it is clear that $[\mJ(t) \cdot \mV(t)]$ is an element of $\mL^2(\Omega)$ for all $t > 0$. 

Next, for each multi-index $\alpha$, $|\alpha| \leq m$, $D^\alpha [\mJ(t) \cdot \mV(t)]$ is defined by
\[
 \langle D^\alpha [\mJ(t) \cdot \mV(t)], \Phi \rangle_{\mL^2(\Omega)}=(-1)^{|\alpha|} \langle [\mJ(t) \cdot \mV(t)], D^\alpha \Phi \rangle_{\mL^2(\Omega)},
\] 
for all $\Phi$ in $\mathbf{C}_0^\infty(\Omega,\R^p)$, the space of infinitely differentiable $\R^p$-valued functions defined on $\Omega$ with compact support. Applying Fubini's theorem again we obtain that the righthand side of the previous equality is equal to
\[
 \langle [D^\alpha \mJ(t) \cdot \mV(t)],\Phi \rangle_{\mL^2(\Omega)}.
\]
This proves that $D^\alpha [\mJ(t) \cdot \mV(t)]$ exists for all $t>0$, all multi-indexes $\alpha$, $|\alpha| \leq m$, and is equal to $[D^\alpha \mJ(t) \cdot \mV(t)]$. To see that it is in $\mL^2(\Omega)$ just apply Fubini's theorem again.
\end{proof}
\section{Fixed points theorems}\label{appendix:LS}
We briefly describe some applications of Brouwer's and Leray-Schauder's degree theories. They are central in the proofs of parts 4 and 5 of proposition \ref{prop:bumps}.

We first recall Schaeffer's theorem and provide a short proof based on the Leray-Schauder degree.
\begin{theorem}[Schaeffer]\label{theo:schaeffer}
 Let X be a real Banach space. Suppose $M:X\rightarrow X$
  is a compact mapping and 
$$\mathcal S = \left\{x\in X\,|\,\exists t\in[0,1]\ such\ that\ x=tM(x)\right\}$$
is bounded. Then $M$ has a fixed point.
\end{theorem}
\begin{proof}
We provide for completeness a short proof based on Leray-Schauder's degree theory.
Taking $r>0$ large enough such that $\mathcal S\subset B^\circ_r$, we define
$m(x,t)=x - t M(x)$ on $\bar B_r\times [0,1]$. Then $0\notin
m(\partial B_r\times [0,1])$ by construction. According to the homotopy invariance of the Leray-Schauder degree
$$ {\rm deg}_{\rm LS}({\rm Id}-M,B_r,0)={\rm deg}_{\rm LS}({\rm Id},B_r,0)=1\neq 0,$$
thus, according to the Kronecker property of the Leray-Schauder degree, there exists a solution to the equation $M(x)=x$.
\end{proof}

\noindent
We can apply this theorem to prove existence of solutions to equation \eqref{eq:bump}.
We consider the function $\tilde{F} : \mF \to \mF$ defined by:
\[
 \tilde{F}(\mV,\lambda)=-F(\mV,\lambda)+\mV=\mJ \cdot \mS(\lambda \mV)+\mI,
\]
where $F$ is defined in equation \eqref{eq:bump}.

It is known that $\tilde{F}$ is a nonlinear compact operator of $\mF$ \cite{faugeras-veltz-etal:09}. 
We can apply Schaeffer's theorem to the function $\tilde{F}$ since it is easy to prove (see \cite{faugeras-veltz-etal:09}) that for all $\mV$ such that $\tilde{F}(\mV,\lambda)=t\mV$ the following holds
\begin{equation}\label{eq:bounded}
 \| \mV \|_\mF \leq t(\lambda \sqrt{p|\Omega|} \|\mJ\|_F+\|\mI\|_\mF),
\end{equation}
where $\norm{\ }_F$ is the Frobenius norm.

Hence for all $\lambda$ there exists $\mV^f_\lambda$ such that
\[
 \tilde{F}(\mV^f_\lambda,\lambda)=\mV^f_\lambda
\]

The following easy consequence of the proof of theorem \ref{theo:schaeffer} is used in the article.
\begin{corollary}\label{coro:bounded}
For each $\lambda \geq 0$ there exists an open bounded set $\mathcal{U}_\lambda$ containing $\mathcal{S}$ (defined in theorem \ref{theo:schaeffer}) such that
${\rm deg}_{\rm LS}({\rm Id} - \tilde{F},\mathcal{U}_\lambda,0)=1$.
\end{corollary}

\noindent
We next give without proof a theorem due to Leray and Schauder.
\begin{theorem}[Leray-Schauder] 
Let $X$ be a real normed space, $J=[a,b]$ and $
  M:X\times J\rightarrow X$ be of the form $Id+m$ with $m:X\times
  J\rightarrow X$ compact on $X \times J$. Let 
$$ \Sigma = \left\{(x,\lambda)\in X\times J\ :
  \ M(x,\lambda)=0\right\}$$
and for each $\lambda\in J$, let 
$$\Sigma_\lambda = \left\{x\in X\ :\ (x,\lambda)\in \Sigma\right\}$$
Assume that $\Sigma_a$ is bounded and
that $${\rm deg}_{\rm LS}(M(.,a),\mathcal{U},0)\neq 0$$
for some open bounded set $\mathcal{U} \supset \Sigma_a$.

Then $\Sigma$ contains a connected component $\mathcal C$ intersecting
$\Sigma_a\times\left\{a\right\}$ and which either intersects
  $\Sigma_b\times\left\{b\right\}$ or is unbounded.
\end{theorem}

\section{Compact operators with simple eigenvalues}\label{appendix:simpleeigen}
\begin{proposition}
\label{prop:DensiteDiag}
For every $m\in\mathbb N^*$, the set of compact operators with $m$ simple first eigenvalues  is dense in the set of compact operators.
\end{proposition}
\begin{proof}
The set $\mathcal{R}_f(\mF)$ of finite dimensional range linear operators is dense in the set $\maH$ of the linear compact operators of $\mF$ \cite{brezis:83}. Thus we only need to prove the theorem for $\mathcal{R}_f(\mF)$. Let us consider $\mJ\in\mathcal{R}_f(\mF)$, and $R(\mJ)=Span(e_1,...,e_k)$, $k \geq m$. Without loss of generality we assume that the first eigenvalue of $\mJ$ has multiplicity two, i.e. $\beta_1=\beta_2$. Its corresponding Jordan block is then 
$$\left[
\begin{array}{cc}
\beta_1 & \varepsilon \\
0  & \beta_1
\end{array}
\right] $$

Then if we define $\mJ_n=\mJ+\frac{1}{n} e_2 \otimes e_2$, we have $\lim_{n \to \infty} \mJ_n = \mJ$. The first two eigenvalues $\beta_1$ and $\beta_1+\frac{1}{n}$ of $\mJ_n$ are simple, i.e. $\mJ_n e_1=\beta_1 e_1$ and $\mJ_n (n\varepsilon e_1+e_2)=(\beta_1+\frac{1}{n})(n\varepsilon e_1+e_2)$. We can do the same for the other eigenvalues, and define an operator $\mJ_n$ with $m$ simple first eigenvalues, a finite dimensional rank, and which is arbitrarily close to $\mJ$.
\end{proof}

\begin{proposition}
 The same proposition holds for operators of the type $\rm{ Id}+\mJ$ where $\mJ$ is a compact operator.
\end{proposition}
\begin{proof}
 It follows from that of \ref{prop:DensiteDiag}.
\end{proof}

\section{Reduction of the activity based model to a finite number of ordinary differential equations}\label{appendix:activity}
We consider the equation for the activity-based model:
\begin{equation}
\label{eq:NMA}
\left\{
\begin{array}{lcl}
\dot\mA&=&-\mL_a \cdot \mA+\mS_\lambda\left(\mJ \cdot \mA+\mI\right) \quad t>0\\
\mA(\cdot,0)&=&\mA_0(\cdot)
\end{array}
\right.
\end{equation}

We recall that $\mL_a\neq \mL$ (see \cite{ermentrout:98}) because they do not have the same biologiocal meaning: One is related to the synaptic time constant and the other to the cell membrane time constant. We let $\mL_a={\rm diag}(\alpha_1,\cdots,\alpha_p)$

We also recall the PG-kernel decomposition of $\mJ=\sum\limits_k X_k\otimes Y_k$ 

Each of the $p$ coordinates $A_i$, $i=1,\cdots,p$ of $\mA$ satisfies 
\[
\dot A_i+\alpha_i A_i=\mS_\lambda(\mJ\cdot\mA +\mI)_i \quad i=1,\cdots,p 
\]
Similarly to the voltage case, let us consider the $p$ finite dimensional subspaces $F_i$, $i=1,\cdots,p$ of $\mG$, where each $F_i$ is generated by the $N$ elements  $Y^i_k,\,k=1,\cdots,N$. We decompose $\mG$ as the direct sum of $F_i$ and its orthogonal complement for $\left\langle ,\right\rangle_2$ $F_i^\bot$, $\mG=F_i \oplus F_i^\bot$ and write $A_i=A_i^\parallel+A_i^\bot$. We note $\prod_i^\parallel$ (respectively $\prod_i^\bot$) the projection from  $\mG$ to $F_i$ (respectively to $F_i^\bot$) parallel to $F_i^\bot$ (respectively to $F_i$). This induces a decomposition of $\mF$ as the direct sum of $F=\prod_{i=1}^p F_i$ and $F^\bot=\prod_{i=1}^p F_i^\bot$ such that 
for each vector $\mA$ of $\mF$ we can write $\mA=\mA^\parallel+\mA^\bot$. By construction we also have
\[
 \mJ \cdot \mA^\bot= 0,
\]
and therefore
 \begin{equation*}
  \begin{cases}
   \dot A_i^\parallel +\alpha_i A_i^\parallel=\prod_i^\parallel \mS_\lambda(\mJ \cdot \mA^\parallel +\mI)_i\\
   \dot A_i^\bot+\alpha_i A_i^\bot=\prod_i^\bot \mS_\lambda(\mJ \cdot \mA_F+\mI)_i
  \end{cases} i=1,\cdots,p
 \end{equation*}
which is a $2p$-dimensional non-autonomous system of ODEs:
\begin{equation*}
  \begin{cases}
   \dot\mA^\parallel+\mL_a \cdot \mA^\parallel=\prod^\parallel \mS_\lambda(\mJ\cdot\mA^\parallel +\mI)\\
   \dot\mA^\bot+\mL_a \cdot \mA^\bot=\prod^\bot\mS_\lambda(\mJ\cdot\mA^\parallel +\mI) 
  \end{cases},
 \end{equation*}
where $\prod^\parallel \mA=(\prod_i^\parallel A_i)_{i=1,\cdots,p}$ and $\prod^\bot \mA=(\prod_i^\bot A_i)_{i=1,\cdots,p}$ are the projections of $\mA$ on $F$ and $F^\bot$.

The first equation is a $p$-dimensional autonomous system of ODEs, which can be solved before solving the second one. 

\subsection{Lemmas}\label{appendix:lemmas}
\begin{lemma}\label{lemma:Sromain}
For all $x,\, \lambda \in \R$ we have
\[
 (S(\lambda x)-S(0))^2 \leq S(\lambda^2 x^2)-S(0)
\]

\end{lemma}
\begin{proof}
We set $X=\lambda x$ and consider two cases.
\begin{description}
 	\item[$X > 1$] We have $e^{-X} > e^{-X^2}$ and therefore $S(X)-1/2 < S(X^2)-1/2$. Moreover, since $S(X)-1/2 < 1$, $(S(X)-1/2)^2 < S(X)-1/2$ and we are done.
	\item[$0< X < 1$] We let $X=\log y$, $1 < y < e$. We therefore have
\[
 S(X)-1/2=\frac{1}{2} \frac{y-1}{y+1} \quad S(X^2)-1/2=\frac{1}{2} \frac{y^{\log y}-1}{y^{\log y}+1}
\]
We consider the expression $(y-1)^2 (y^{\log y}+1)-2(y+1)^2(y^{\log y}-1)$ and prove it is negative. Because $y^{\log y} < y$ it is upperbounded by $(y-1)^2 (y^{\log y}+1)-2(y^{\log y}+1)^2(y^{\log y}-1)$ which
has the sign of $(y-1)^2-2(y^{2\log y}-1)$. The last expression is upperbounded by $(y-1)^2-2(y^{\log y}-1)=(e^X-1)^2-2(e^{X^2}-1)$ which is negative for $0 < X <1$.
\end{description}
\end{proof}

\begin{proposition}\label{prop:ineq1}
The solutions of equation \eqref{eq:bump} satisfy the following inequalities for all $\lambda \geq 0$
\[
 \norm{\mV^f_\lambda}_\mF \leq \sqrt{p|\Omega|} \left\|\mJ\right\|_\mF+\norm{\mI_{\rm ext}}_\mF \overset{\rm def}{=}B_1
\]
\[
 \left\|\mV^f_\lambda-\mV^f_0\right\|_\mF\leq  \frac{1}{2} \sqrt{p|\Omega|} \left\|\mJ\right\|_\mF \overset{\rm def}{=}B_2,
\]
as well as
\[
\left\|\mV^f_\lambda-\mV^f_0\right\|_\mF\leq \frac{\lambda}{4} \left\|\mJ\right\|_\mF B_1
\]
\end{proposition}
\begin{proof}
The first inequality is a straightforward consequence of equation \eqref{eq:bump}, taking the $\mF$-norm and using the fact that $0 \leq S(x) \leq 1$ for all $x \in \R$.http://fidji.inria.fr/biblio
For the second one we write $\mV^f_\lambda-\mV^f_0=\mJ \cdot (\mS(\lambda \mV^f_\lambda)-\mS(0))$, take the $\mF$-norm of both sides of the equality and use the Cauchy-Schwarz inequality. We find $\norm{\mV^f_\lambda-\mV^f_0}_\mF \leq \norm{\mJ}_\mF \cdot \norm{\mS(\lambda \mV^f_\lambda)-\mS(0)}_{{\rm L}^2(\Omega,\mathbb R^p)}$. But since $\forall x\in \R,\ -\frac{1}{2} \leq S(x)-S(0)\leq \frac{1}{2}$, we have $\norm{\mS(\mV^f_\lambda)-\mS(0)}_{{\rm L}^2(\Omega,\mathbb R^p)} \leq  \frac{1}{2}\sqrt{p|\Omega|}$, which proves the first inequality.

The third inequality can be obtained as follows. It is easy to see that $S(\lambda x)-S(0) \leq \frac{\lambda}{4} |x|$ for all $x \in \R$ and all $\lambda \geq 0$. This implies that $\norm{\mS(\mV^f_\lambda)-\mS(0)} \leq \frac{\lambda}{4} \norm{\mV^f_\lambda}$. The first inequality yields the third.
\end{proof}

\section{An equality for the adjoint operator}\label{appendix:ajoint}
We prove proposition \ref{prop:inner}.
\begin{proof}
 We prove the proposition in the case $m=p=1$. By definition of $\mJ^*_\mF$ we have
\[
 \left\langle \mU, \mJ^*_\mF \cdot \mV \right\rangle_\mF=\left\langle \mJ \cdot \mU, \mV \right\rangle_\mF \quad \forall \mU,\,\mV \in \mF
\]
We note $\mX=\mJ^*_\mF \cdot \mV$ and rewrite the previous equation as
\[
 \left\langle \mU, \mX \right\rangle_2+\left\langle \nabla \mU, \nabla \mX  \right\rangle_2=\left\langle \mJ \cdot \mU, \mV \right\rangle_\mF
\]
We next rewrite the righthand side
\[
 \left\langle \mJ \cdot \mU, \mV \right\rangle_\mF=\left\langle \mJ \cdot \mU, \mV \right\rangle_2+\left\langle \nabla \left(\mJ \cdot \mU \right), \nabla \mV \right\rangle_2
\]
The crucial step is to observe that
\[
 \left\langle \nabla \left(\mJ \cdot \mU \right), \nabla \mV \right\rangle_2=\left\langle \mU, (\nabla_1 \mJ)^*_2 \cdot \nabla \mV \right\rangle_2,
\]
where $\nabla_1$ indicates the derivative with respect to the first variable, i.e.
\[
 (\nabla_1 \mJ)^*_2 \cdot \nabla \mV (\mr)=\sum_{i=1}^d \int_\Omega \frac{\partial \mJ}{\partial r_i'} (\mr',\mr) \frac{\partial \mV}{\partial r_i'}(\mr')\,d\mr'
\]
Let us fix $\mV$ in $\mF$. $\mX$ is the solution of
\[
 \left\langle \mU, \mX \right\rangle_2+\left\langle \nabla \mU, \nabla \mX  \right\rangle_2= \left\langle \mU, f \right\rangle_2 \quad \forall \mU \in \mF
\]
Hence $\mX$ is the unique weak solution in $\mF$ of the homogeneous Neumann problem
\[
 \left\{
\begin{array}{ccc}
 -\Delta \mX+\mX & = & f \quad \text{on} \quad \Omega\\
\frac{\partial \mX}{\partial \nu} & = & 0 \quad \text{on} \quad \partial \Omega
\end{array}
\right.
\]
$f$ is the element of $\mL^2(\Omega)$ defined by
\[
 f=\mJ^*_2 \cdot \mV+ (\nabla_1 \mJ)^*_2 \cdot \nabla \mV+
\]
The regularity properties of the solutions of elliptic equations \cite[Chapter 6]{evans:98} imply that $\mX$ is in $\mW^{2,2}(\Omega)$. Let us denote by $\mathcal{D}$ the differential operator, defined in $\mW^{2,2}(\Omega)$ by
\[
 \mathcal{D}=-\Delta+{\rm Id}
\]
We now choose $\mV=e^*_\mF$. This implies that $\mX=\lambda e^*_\mF$ and hence $\frac{\partial e^*_\mF}{\partial \nu}=0$ in $\partial \Omega$. Furthermore we have
\[
 \lambda \mathcal{D} \cdot e^*_\mF=\mJ^*_2 \cdot e^*_\mF+(\nabla_1 \mJ)^*_2 \cdot \nabla e^*_\mF
\]
Green's formula with the condition $\frac{\partial e^*_\mF}{\partial \nu}=0$ show that
\[
 (\nabla_1 \mJ)^*_2 \cdot \nabla e^*_\mF=-\mJ^*_2 \cdot \Delta e^*_\mF,
\]
and therefore we have
\[
 \lambda \mathcal{D} \cdot e^*_\mF= \mJ^*_2 \cdot (\mathcal{D} \cdot e^*_\mF),
\]
which shows that $\mathcal{D} \cdot e^*_\mF$ is an eigenvector of $\mJ^*_2$ associated with the eigenvalue $\lambda$. We can choose $e^*_{L^2}$ and $e^*_\mF$ so that
\[
 \mathcal{D} \cdot e^*_\mF=e^*_{L^2}.
\]
From Green's formula
\[
 \left\langle \mU, e^*_\mF \right\rangle_\mF=\left\langle \mU, \mathcal{D} \cdot e^*_\mF \right\rangle_2-\int_{\partial \Omega} \frac{\partial e^*_\mF}{\partial \nu} \mU \, d\sigma,
\]
and the conclusion follows from the fact that $\mathcal{D} \cdot e^*_\mF=e^*_{L^2}$ and $\frac{\partial e^*_\mF}{\partial \nu}=0$ on $\partial \Omega$.
\end{proof}

\newpage
\bibliographystyle{siam}

\bibliography{/local/home/faugeras/latex/string,/local/home/faugeras/latex/odyssee}

\begin{thebibliography}{10}

\bibitem{adams:75}
{\sc R.~Adams}, {\em Sobolev spaces}, vol.~65 of Pure and Applied Mathematics,
  Series of Monographs and Textbooks, Academic Press, Inc., New York, San
  Francisco, London, 1975.

\bibitem{amari:77}
{\sc S.-I. Amari}, {\em Dynamics of pattern formation in lateral-inhibition
  type neural fields}, Biological Cybernetics, 27 (1977), pp.~77--87.

\bibitem{atay-hutt:05}
{\sc Fatihcan~M. Atay and Axel Hutt}, {\em Stability and bifurcations in neural
  fields with finite propagation speed and general connectivity}, SIAM Journal
  on Applied Mathematics, 65 (2005), pp.~644--666.

\bibitem{ben-yishai-bar-or-etal:95}
{\sc R.~Ben-Yishai, RL~Bar-Or, and H.~Sompolinsky}, {\em Theory of orientation
  tuning in visual cortex}, Proceedings of the National Academy of Sciences, 92
  (1995), pp.~3844--3848.

\bibitem{blomquist-wyller-etal:05}
{\sc P.~Blomquist, J.~Wyller, and G.T. Einevoll}, {\em Localized activity
  patterns in two-population neuronal networks}, Physica D, 206 (2005),
  pp.~180--212.

\bibitem{bressloff:05}
{\sc P.~Bressloff}, {\em Spontaneous symmetry breaking in self--organizing
  neural fields}, Biological Cybernetics, 93 (2005), pp.~256--274.

\bibitem{bressloff-bressloff-etal:00}
{\sc PC~Bressloff, NW~Bressloff, and JD~Cowan}, {\em Dynamical mechanism for
  sharp orientation tuning in an integrate-and-fire model of a cortical
  hypercolumn}, Neural computation, 12 (2000), pp.~2473--2511.

\bibitem{bressloff-cowan:02c}
{\sc P.C. Bressloff and J.D. Cowan}, {\em An amplitude equation approach to
  contextual effects in visual cortex}, Neural computation, 14 (2002),
  pp.~493--525.

\bibitem{bressloff-cowan-etal:01}
{\sc P.C. Bressloff, J.D. Cowan, M.~Golubitsky, P.J. Thomas, and M.C. Wiener},
  {\em Geometric visual hallucinations, {E}uclidean symmetry and the functional
  architecture of striate cortex}, Phil. Trans. R. Soc. Lond. B, 306 (2001),
  pp.~299--330.

\bibitem{bressloff-folias-etal:03}
{\sc PC~Bressloff, SE~Folias, A.~Prat, and Y.X. Li}, {\em Oscillatory waves in
  inhomogeneous neural media}, Physical Review Letters, 91 (2003), p.~178101.

\bibitem{brezis:83}
{\sc H.~Brezis}, {\em Analyse fonctionnelle. Th\'eorie et applications},
  Masson, 1983.

\bibitem{choquet:69}
{\sc G.~Choquet}, {\em Cours d'Analyse}, vol.~II, Masson, 1969.

\bibitem{colby-duhamel-etal:95}
{\sc C.L. Colby, J.R. Duhamel, and M.E. Goldberg}, {\em Oculocentric spatial
  representation in parietal cortex}, Cereb. Cortex, 5 (1995), pp.~470--481.

\bibitem{coombes:05}
{\sc Stephen Coombes}, {\em Waves, bumps, and patterns in neural fields
  theories}, Biological Cybernetics, 93 (2005), pp.~91--108.

\bibitem{coombes-owen:04}
{\sc S.~Coombes and M.~R. Owen}, {\em Evans functions for integral neural field
  equations with heaviside firing rate function}, SIAM Journal on Applied
  Dynamical Systems, 3 (2004), pp.~574--600.

\bibitem{coombes-owen:05}
\leavevmode\vrule height 2pt depth -1.6pt width 23pt, {\em Bumps, breathers,
  and waves in a neural network with spike frequency adaptation}, Phys. Rev.
  Lett., 94 (2005).

\bibitem{dayan-abbott:01}
{\sc P.~Dayan and L.~F. Abbott}, {\em Theoretical Neuroscience : Computational
  and Mathematical Modeling of Neural Systems}, MIT Press, 2001.

\bibitem{ermentrout:98}
{\sc Bard Ermentrout}, {\em Neural networks as spatio-temporal pattern-forming
  systems}, Reports on Progress in Physics, 61 (1998), pp.~353--430.

\bibitem{ermentrout-cowan:80}
{\sc GB~Ermentrout and JD~Cowan}, {\em Large scale spatially organized activity
  in neural nets}, SIAM Journal on Applied Mathematics,  (1980), pp.~1--21.

\bibitem{evans:98}
{\sc L.C. Evans}, {\em {P}artial {D}ifferential {E}quations}, vol.~19 of
  Graduate Studies in Mathematics, Proceedings of the American Mathematical
  Society, 1998.

\bibitem{faugeras-grimbert-etal:08}
{\sc O.~Faugeras, F.~Grimbert, and J.-J. Slotine}, {\em Abolute stability and
  complete synchronization in a class of neural fields models}, SIAM Journal of
  Applied Mathematics, 61 (2008), pp.~205--250.

\bibitem{faugeras-touboul-etal:09}
{\sc O.~Faugeras, J.~Touboul, and B.~Cessac}, {\em A constructive mean field
  analysis of multi population neural networks with random synaptic weights and
  stochastic inputs}, Frontiers in Computational Neuroscience, 3 (2009).

\bibitem{faugeras-veltz-etal:09}
{\sc O.~Faugeras, R.~Veltz, and F.~Grimbert}, {\em Persistent neural states:
  stationary localized activity patterns in nonlinear continuous n-population,
  q-dimensional neural networks}, Neural Computation, 21 (2009), pp.~147--187.

\bibitem{folias-bressloff:05}
{\sc SE~Folias and PC~Bressloff}, {\em Breathers in two-dimensional neural
  media}, Physical Review Letters, 95 (2005), p.~208107.

\bibitem{funahashi-bruce-etal:89}
{\sc S.~Funahashi, C.J. Bruce, and P.S. Goldman-Rakic}, {\em Mnemonic coding of
  visual space in the monkey's dorsolateral prefrontal cortex}, J.
  Neurophysiol., 61 (1989), pp.~331--349.

\bibitem{guckenheimer-holmes:83}
{\sc J.~Guckenheimer and P.~J. Holmes}, {\em Nonlinear Oscillations, Dynamical
  Systems and Bifurcations of Vector Fields}, vol.~42 of Applied mathematical
  sciences, Springer, 1983.

\bibitem{hansel-sompolinsky:97}
{\sc D.~Hansel and H.~Sompolinsky}, {\em Modeling feature selectivity in local
  cortical circuits}, Methods of neuronal modeling,  (1997), pp.~499--567.

\bibitem{haragus-iooss:09}
{\sc M.~Haragus and G.~Iooss}, {\em Local bifurcations, center manifolds, and
  normal forms in infinite dimensional systems}, EDP Sci., 2009.
\newblock To appear.

\bibitem{henry:81}
{\sc D.~Henry}, {\em Geometric Theory of Semilinear Parabolic Equations},
  vol.~840 of LNM, Springer-Verlag, 1981.

\bibitem{kielhofer:03}
{\sc H.Kielh{\"o}fer}, {\em Bifurcation Theory: An Introduction with
  Applications to PDEs}, Springer, 2003.

\bibitem{kato:95}
{\sc T.~Kato}, {\em {Perturbation Theory for Linear Operators}}, Springer,
  1995.

\bibitem{kuznetsov:98}
{\sc Yuri~A. Kuznetsov}, {\em Elements of Applied Bifurcation Theory}, Applied
  Mathematical Sciences, Springer, 2nd~ed., 1998.

\bibitem{laing-troy-etal:02}
{\sc C.L. Laing, W.C. Troy, B.~Gutkin, and G.B. Ermentrout}, {\em Multiple
  bumps in a neuronal model of working memory}, SIAM J. Appl. Math., 63 (2002),
  pp.~62--97.

\bibitem{laing-troy:03}
{\sc Carlo~R. Laing and William~C. Troy}, {\em {PDE} methods for nonlocal
  models}, SIAM Journal on Applied Dynamical Systems, 2 (2003), pp.~487--516.

\bibitem{ma-wang:05}
{\sc Tian Ma and Shouhong Wang}, {\em Bifurcarion theory and applications},
  vol.~53 of Nonlinear Science, World Scientific, 2005.

\bibitem{miller-erickson-etal:96}
{\sc E.K. Miller, C.A. Erickson, and R.~Desimone}, {\em Neural mechanisms of
  visual working memory in prefrontal cortex of the {M}acaque}, J. Neurosci.,
  16 (1996), pp.~5154--5167.

\bibitem{nishiura-mimura:89}
{\sc Yasumasa Nishiura and Masayasu Mimura}, {\em Layer oscillations in
  reaction-diffusion systems}, SIAM Journal on Applied Mathematics, 49 (1989),
  pp.~481--514.

\bibitem{rabinowitz:71}
{\sc P.H. Rabinowitz}, {\em Some global results for nonlinear eigenvalue
  problems}, J. Funct. Anal, 7 (1971), pp.~487--513.

\bibitem{sala-heroux-etal:04}
{\sc Marzio Sala, Michael~A. Heroux, and David~M. Day}, {\em Trilinos
  tutorial}, Tech. Report SAND2004-2189, Sandia National Laboratories, 2004.

\bibitem{shriki-hansel-etal:03}
{\sc O.~Shriki, D.~Hansel, and H.~Sompolinsky}, {\em Rate models for
  conductance-based cortical neuronal networks}, Neural Computation, 15 (2003),
  pp.~1809--1841.

\bibitem{tricomi:85}
{\sc F.G. Tricomi}, {\em Integral Equations}, Dover, 1985.
\newblock Reprint.

\bibitem{venkov-coombes-etal:07}
{\sc N.A. Venkov, S.~Coombes, and P.C. Matthews}, {\em Dynamic instabilities in
  scalar neural field equations with space-dependent delays}, Physica D:
  Nonlinear Phenomena, 232 (2007), pp.~1--15.

\bibitem{wilson-cowan:73}
{\sc H.R. Wilson and J.D. Cowan}, {\em A mathematical theory of the functional
  dynamics of cortical and thalamic nervous tissue}, Biological Cybernetics, 13
  (1973), pp.~55--80.

\bibitem{yosida:95}
{\sc K.~Yosida}, {\em Functional analysis. reprint of the sixth (1980) edition.
  classics in mathematics}, 1995.

\end{thebibliography}

\end{document}